\documentclass[11pt]{article}
\usepackage{graphicx, color}
\usepackage{amsmath,amsfonts,amssymb}
\usepackage{graphicx}

\def\ds{\displaystyle}
\def\forall{\hbox{for all}~}
\def\L{{\bf L}}
\def\H{{\cal H}}

\def\bfn{{\bf n}}

\def\ve{\varepsilon}

\def\avint{-\!\!\!\!\!\!\int}

\def\E{{\cal E}}
\def\I{{\cal I}}

\def\S{{\cal S}}

\def\wto{\rightharpoonup}
\def\R{{\mathbb R}}

\def\implies{\Longrightarrow}
\def\vp{\varphi}

\def\vs{\vskip 2em}

\def\v{\vskip 1em}

\def\begi{\begin{itemize}}
\def\endi{\end{itemize}}
\def\C{{\cal C}}
\def\M{{\cal M}}
\def\ov{\overline}
\def\Tilde{\widetilde}
\def\Hat{\widehat}
\def\bega{\begin{array}}
\def\enda{\end{array}}
\def\meas{\hbox{meas}}

\def\bel{\begin{equation}\label}
\def\eeq{\end{equation}}
\def\sqr#1#2{\vbox{\hrule height .#2pt
\hbox{\vrule width .#2pt height #1pt \kern #1pt
\vrule width .#2pt}\hrule height .#2pt }}
\def\square{\sqr74}
\def\endproof{\hphantom{MM}\hfill\llap{$\square$}\goodbreak}

\newtheorem{theorem}{Theorem}[section]
\newtheorem{lemma}[theorem]{Lemma}
\newtheorem{proposition}[theorem]{Proposition}
\newtheorem{corollary}[theorem]{Corollary}
\newtheorem{definition}[theorem]{Definition}
\newtheorem{remark}[theorem]{Remark}

\begin{document}
\title{\bf  Variational Problems for Tree Roots and Branches}
\vs
\author{Alberto Bressan, Michele Palladino, and Qing Sun
\\
\,
\\
Department of Mathematics, Penn State University, \\
University Park, Pa.~16802, USA.\\
\,
\\
e-mails: axb62@psu.edu, ~mup26@psu.edu,~qxs15@psu.edu}
\maketitle
\begin{abstract} 
 This paper studies two classes of variational problems 
introduced in \cite{BSu}, related to the
optimal shapes of tree roots and branches.
Given a measure $\mu$ describing the distribution of leaves, 
a  {\it sunlight functional} $\S(\mu)$ computes the total 
amount of light captured by the leaves.
For a measure $\mu$ describing the distribution of root hair cells, 
a {\it harvest functional} $\H(\mu)$ computes the total amount of water and nutrients
gathered by the roots.  
In both cases, we seek a measure $\mu$ that maximizes 
these functionals subject to a 
ramified transportation cost, for transporting 
nutrients from the roots to the trunk or from the trunk to the leaves. 

Compared with \cite{BSu}, here we do not impose any a priori bound 
on the total mass of the optimal measure $\mu$, and 
more careful a priori estimates are thus required.
In the unconstrained optimization problem for branches, we prove that 
an optimal measure exists, with bounded support and bounded total mass.
In the unconstrained problem for tree roots, we prove that an optimal measure exists,
with bounded support but possibly unbounded total mass.
The last section of the paper analyzes how the size of the optimal tree
depends on the parameters defining the various functionals.
\end{abstract}
\v
         

\vs
\section{Introduction}
\label{s:0}
\setcounter{equation}{0}
In the recent paper \cite{BSu}, two of the authors introduced a family of variational problems,
aimed at characterizing optimal shapes of tree roots and branches.
All these optimization  problems take place in a space of positive measures on a $d$-dimensional
space $\R^d$.  
In the case of roots, calling $\mu$ the distribution of root hair cells,
one seeks to maximize the total amount of water and nutrients harvested by the roots, 
minus a cost for transporting these nutrients to the base of the trunk. 
In the case of branches,
calling $\mu$ the distribution of leaves, 
one seeks to maximize the total sunlight captured
by the leaves, minus a cost for transporting water and nutrients from the base of the trunk 
to the tip of every branch.   

The main results in \cite{BSu} established the semicontinuity of the relevant functionals
and the existence of optimal solutions, under a constraint on the total mass 
of the measure $\mu$.  In essence, by fixing the total mass $\mu(\R^d)$ 
one  prescribes the {\it size} of the tree.   In turn, the maximization problem determines
an {\it optimal shape}.  

In the present paper we study the corresponding {\it unconstrained} optimization
problems, without any a priori bound on the total mass of the measure $\mu$.
Roughly speaking, this aims at
determining the {\it optimal size} of a tree, in addition to its optimal shape.

Compared with \cite{BSu}, proving the existence 
of optimal solutions for the unconstrained problems requires a much more careful analysis.
Following the direct method of the Calculus of Variations, we 
consider a maximizing sequence of measures $(\mu_k)_{k\geq 1}$.
Two main issues arise.
\begi
\item[(i)] First, one needs to establish an a priori bound on the support 
of the measures $\mu_k$.   At first sight this looks easy, because if a measure
contains some mass far away from the origin, its transportation cost will be very large.
However, since we are here considering a {\it ramified transportation cost} 
\cite{BCM,  MMS, X03, X4}, there is an {\it economy of scale}:  
as the total transported mass increases without bound, the unit cost decreases to zero.
For this reason, in order to achieve a uniform bound on  Supp$(\mu_k)$, 
we first establish an a priori bound on the transportation cost. 
At a second stage, this yields a bound on the total payoff.  Finally, we obtain
an estimate of the support of the optimal measure.

\item[(ii)] Next, we seek an a priori bound on  the total mass $\mu_k(\R^d)$.  
This does not follow from a bound on the transportation cost, because 
as $k\to\infty$  the measures $\mu_k$ may concentrate more and more mass
in a small neighborhood of the origin. Concerning  the optimization problem for branches,
our analysis yields the existence of an optimal measure $\mu$ such that 
$\mu(\R^d)<+\infty$.  On the other hand,  
in  the optimization problem for tree roots, 
we prove that an optimal measure $\mu$ exists, with bounded support 
but possibly 
unbounded total mass.  Indeed,  for any $\rho>0$ we can show that
$\mu\bigl(\{x\in\R^d\,;~~|x|>\rho\}\bigr)<+\infty$.   However, we cannot rule out 
the possibility that $\mu(\R^d\setminus \{0\})=+\infty$.
\endi

The remainder of the paper is organized as follows. 
Section~2 reviews the three main ingredients of
our variational problems: the {\it sunlight functional}, the {\it harvest functional}, 
and the {\it ramified transportation cost}.
In Section~3 we prove the existence of a bounded measure $\mu$ which solves the 
unconstrained optimization problem for tree branches.   The proof relies on the construction of a
maximizing sequence of measures $(\mu_k)_{k\geq 1}$ with uniformly bounded support
and uniformly bounded total mass. In this direction, a key step is to prove a uniform bound
on the ramified transportation cost for all measures $\mu_k$.
Section~4 deals with the unconstrained optimization problem for tree roots.
The existence of an optimal measure $\mu$  is proved, with bounded support but possibly infinite total  mass.
Finally, in Section~5 we discuss how the optimal size of tree roots and branches
is affected by the various parameters appearing in the equations.  Here the key step is 
to analyze  how  the various functionals behave under a rescaling of coordinates. 

The theory of ramified transport for general measures was developed independently in 
\cite{MMS} and   \cite{X03}.  See also \cite{BCM} for a comprehensive introduction, and
\cite{X4} for a survey of the field.  Further results on optimal 
ramified transport can be found in \cite{BCM1, BS, MS, MoS, S1}.    An interesting 
computational
approach, based on Gamma-convergence, has been developed in \cite{OS,  S2}.
A geometric optimization problem involving a ramified transportation cost was recently
studied in \cite{PSX}.
The ``sunlight functional" was introduced in \cite{BSu}, in a slightly more general setting
which also takes into account the presence of external vegetation.    
The ``harvest functional", in a space of Radon measures, 
was first studied in \cite{BCS} in connection with a problem of optimal harvesting of marine resources.
\v

\section{Review of the basic functionals}
\label{s:1}
\setcounter{equation}{1}
Given a positive, bounded Radon measure $\mu$ on $\R^d$,
 three functionals were considered in \cite{BSu}.  The  corresponding optimization problems determine
 the optimal configurations of roots and branches of a tree.

\centerline{\hbox{\includegraphics[width=9cm]{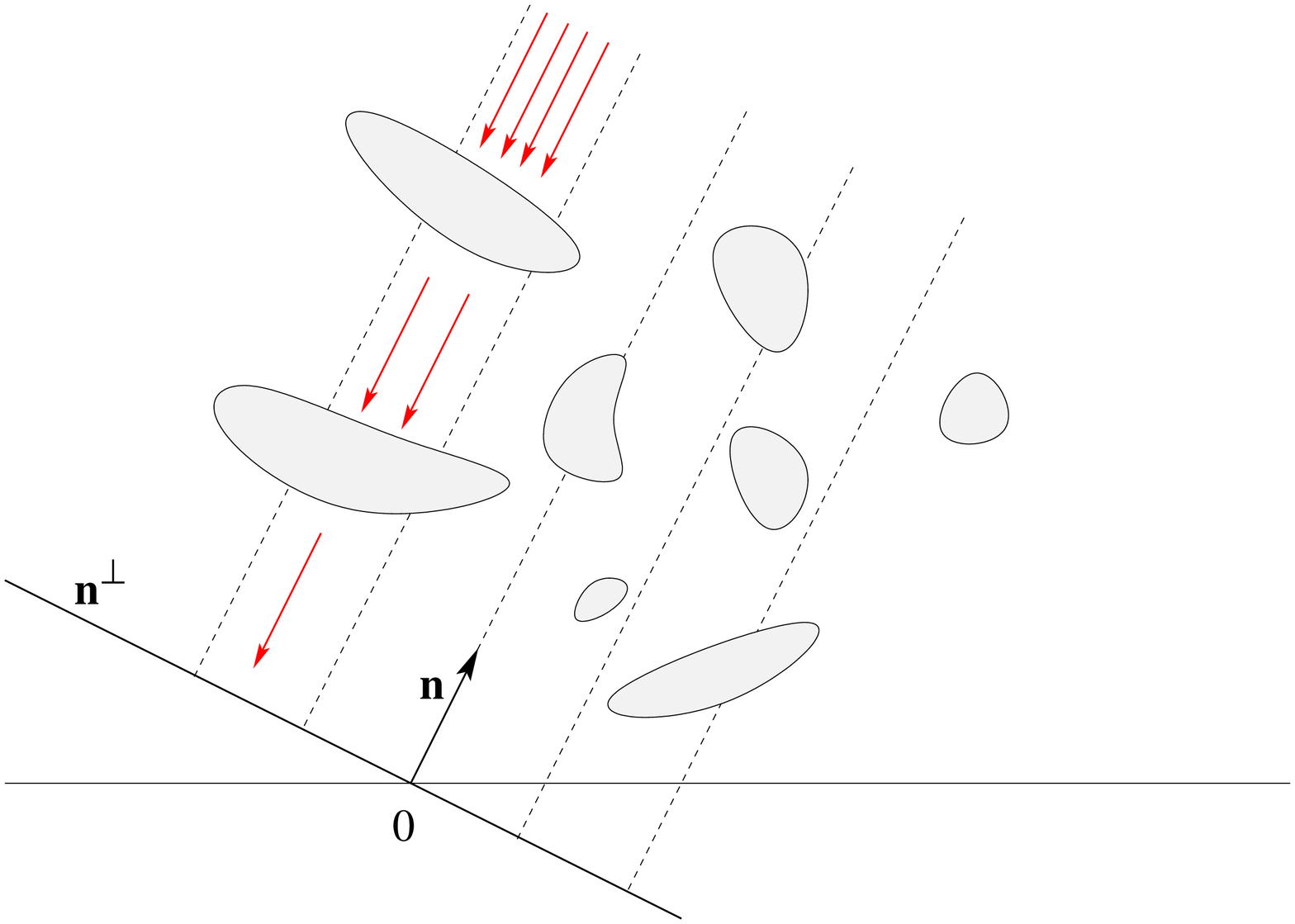}}}
\begin{figure}[ht]
\caption{\small  Sunlight arrives from the direction 
parallel to $\bfn$.  
Part of it is absorbed by the measure $\mu$, 
supported on the grey regions.}
\label{f:ir10}
\end{figure}

\subsection{A sunlight functional}
Let $\mu$ be a positive, bounded Radon measure on $\R^d$.   
Thinking of $\mu$ as the density of leaves on a tree, 
we seek a functional $\S(\mu)$ describing 
the total amount of sunlight absorbed by the leaves.
As shown  in Fig.~\ref{f:ir10}, fix a unit vector
$$\bfn~\in~S^{d-1}~=~\{ x\in \R^d\,;~~|x|=1\}$$ 
and assume that all  
light rays come parallel to $\bfn$.  
Call $E_\bfn^\perp$ the $(d-1)$-dimensional subspace perpendicular to $\bfn$ 
and let $\pi_\bfn:\R^d\mapsto E_\bfn^\perp$ be the perpendicular projection.   Each point $x\in \R^d$ can thus be expressed 
uniquely as
\bel{perp}
x~=~y + s\bfn\eeq
with $y\in  E_\bfn^\perp$ and $s\in\R$.

On the perpendicular subspace $E_\bfn^\perp$ consider the projected measure $\mu^\bfn$, defined by setting
\bel{mupro}\mu^\bfn(A)~=~\mu\Big(\bigl\{ x\in\R^d\,;~~\pi_\bfn(x)\in A\bigr\}\Big).\eeq
Call $\Phi^\bfn$ the density of the absolutely continuous part of $\mu^\bfn$
w.r.t.~the $(d-1)$-dimensional Lebesgue measure on $E_\bfn^\perp$.
\v
\begin{definition}
The total amount of sunshine from the direction $\bfn$ captured by a measure 
$\mu$ on $\R^d$ is defined as
\bel{SSn}
\S^\bfn(\mu)~\doteq~
\int_{E_\bfn^\perp}\Big(1- \exp\bigl\{ - \Phi^\bfn(y)\bigr\}\Big)
\, dy\,.\eeq
Given an integrable function $\eta\in \L^1(\S^{d-1})$, 
the total sunshine absorbed by $\mu$ from all directions
is defined as
\bel{SS2}
\S^\eta(\mu)~\doteq~
\int_{S^{d-1}}\left(\int_{E_\bfn^\perp}\Big(1- \exp\bigl\{ - \Phi^\bfn(y)\bigr\}\Big)
\, dy\right) \eta(\bfn)\,d\bfn\,.\eeq
\end{definition}
\v
We think of  $\eta(\bfn)$ as the intensity of light
coming from  the direction $\bfn$.  We recall two estimates proved in \cite{BSu}.

\begin{lemma}
For any  positive Radon measure $\mu$ on $\R^d$, one has
	\bel{Sbo1}\S^\eta(\mu)~\leq ~\|\eta\|_{\L^1}\cdot \mu(\R^d).\eeq
	If $\mu$ is supported inside a closed ball with radius $r$,  calling
	$\omega_{d-1}$ the surface of the unit sphere in $\R^d$, one has
	\bel{Sbo2} \S^\eta(\mu)~\leq ~\|\eta\|_{\L^1}\cdot \omega_{d-1}\, r^{d-1}.\eeq
\end{lemma}


\v
\subsection{Harvest functionals}
We now consider a utility functional 
associated with roots.  Here the main goal is to collect
moisture and nutrients from the ground. 
To model the efficiency of a root, in the following we let $u(x)$ be the density 
of water+nutrients at the point $x$, and consider a positive Radon measure $\mu$ 
describing the distribution  of root hair. 
\v
Consider the half space 
$\Omega\doteq \{x=(x_1,\ldots, x_d)\,;~x_d<0\}$.
 Let $\mu$ be a positive, bounded Radon measure supported 
 on the closure $\ov\Omega$, such that 
$\mu(V)=0$ for every set $V$ having zero capacity.
Consider the elliptic problem with measure source
\bel{emu} \Delta u + f(u) - u\,\mu~=~0\eeq
and 
 Neumann boundary conditions
\bel{NBC}
 \partial_{\bfn(x)}  u~=~0\qquad \hbox{on}~~\partial\Omega\,.\eeq
 By $\bfn(x)$ we denote the unit outer normal vector
at the boundary point $x\in\partial\Omega$, while
$\partial_\bfn u$ is the derivative of $u$ in the
normal direction.   Of course, in this case (\ref{NBC}) simply means
$$ x_d\,=\,0\qquad\implies\qquad {\partial\over\partial x_d} \,  u~=~0.$$

   If $\mu$ is a general measure
and $u$ is a discontinuous function,  the integral (\ref{HF1}) may not be well defined.
To resolve this issue,  calling
 $$\avint_V u\, dx~=~{1\over \hbox{meas}(V)}\,\int_V u\, dx$$
the average value of $u$ on a set $V$, for each $x\in \ov\Omega$ we consider the limit
\bel{upv}
u(x)~=~\lim_{r\downarrow 0}~\avint_{\Omega\cap B(x,r)} u(y)\, dy.\eeq
As proved in  \cite{FZ}, if $u\in H^1(\Omega)$ then
the above limit exists at all points $x\in \ov\Omega$ with the possible exception of a set 
whose capacity is zero. If the measure $\mu$ satisfies (A3), the integral (\ref{HF1}) is thus well defined.
Our present setting is actually even better, because in (\ref{emu}) $u$ and $\mu$ are positive while $f$ is bounded.
Therefore, if the constant $C$ is chosen large enough, the function $u+C|x|^2$ is subharmonic.  As a consequence, 
the limit (\ref{upv}) is well defined at every point $x\in \ov\Omega$.
 \v
Elliptic problems with measure data have been 
studied in several papers \cite{BG1, BGO, DMOP}
and are now fairly well understood. 
A key fact is that, roughly speaking, the Laplace operator ``does not see" sets with zero capacity.
Following \cite{BG1, BGO} we thus call $\M_b$ the set of all bounded Radon measures on $\ov\Omega$.
Moreover, we denote by $\M_0\subset\M_b$ the family of measures which vanish on Borel sets with zero capacity,
so that
\bel{mu0}\hbox{cap}_2(V)~=~0\qquad\implies\qquad \mu(V)~=~0.\eeq
For the definition and basic properties of capacity we refer to \cite{EG}.
Every measure $\mu\in\M_b$ can be uniquely decomposed as a sum
\bel{mudec}\mu~=~\mu_0 + \mu_s,\eeq
where $\mu_0\in \M_0$ while the measure $\mu_s$ is supported on a set with zero capacity.
In the definition of solutions, the presence of the singular measure $\mu_s$ is disregarded.
  \v
 \begin{definition} Let $\mu$ be a measure in $\M_b$,  decomposed as in (\ref{mudec}).   
 A function
$ u\in\L^\infty(\Omega)\cap H^1(\Omega)$, with pointwise values given by (\ref{upv}), is a solution to  
the elliptic problem (\ref{emu})-(\ref{NBC}) if
 \bel{nsol}
 -\int_\Omega \nabla u\cdot \nabla\vp\, dx+\int_\Omega f(u)\vp\, dx -\int_{\ov\Omega} u\vp\, d\mu_0~=~0\eeq
 for every test function $\vp\in \C^\infty_c(\R^d)$.
\end{definition}

In connection with a solution $u$ of (\ref{emu}),
the total harvest  is defined as
\bel{HF1}\H(u,\mu)~\doteq~\int_{\ov\Omega}~ u\, d\mu\,.\eeq
Throughout the following we
 assume
\begi
\item[{\bf (A1)}]  $f:[0,M]\mapsto\R$ is a $\C^2$ function such that, for 
some constants $M,K$,
\bel{fp1} f(M)=0,\qquad \quad 
0\leq
f(u)\leq K,\qquad f''(u)< 0 \quad\forall u\in [0,M].\eeq
\endi
\v
\subsection{Optimal irrigation plans}

Given $\alpha\in [0,1]$ and a positive measure $\mu$ on $\R^d$, the minimum
cost for $\alpha$-irrigating the measure $\mu$ from the  
origin will be denoted by $\I^\alpha(\mu)$.
Following Maddalena, Morel, and Solimini \cite{MMS}, this 
 can be described as follows.
Let $M=\mu(\R^d)$ be the total mass to be transported and let 
$\Theta=[0,M]$.
We think of each $\theta\in \Theta$ as a ``water particle".
A measurable map 
\bel{iplan}
\chi:\Theta\times \R_+~\mapsto~ \R^d\eeq
is called an {\it admissible irrigation plan}
if 
\begi
\item[(i)] For every $\theta\in \Theta$, the map
$t\mapsto \chi(\theta,t)$ is Lipschitz continuous. 
More precisely, for each $\theta$ there exists a stopping time $T(\theta)$ such that, calling 
$$\dot \chi(\theta,t)~=~{\partial\over\partial t} ~\chi(\theta,t)$$
the partial derivative w.r.t.~time, one has
\bel{stime}\bigl|\dot \chi(\theta,t)\bigr|~=~\left\{ \bega{rl} 1\qquad &\quad\hbox{for a.e.}
~t\in [0, T(\theta)],\\[3mm]
0\qquad &\quad\hbox{for}
~t\geq T(\theta).\enda\right.\eeq
\item[(ii)] At time $t=0$ all particles are at the origin:
$\chi(\theta,0)={\bf 0}$ for all $\theta\in\Theta$.
\item[(iii)] The push-forward of the Lebesgue measure on $[0,M]$ through the map $\theta\mapsto 
\chi(\theta,T(\theta))$ coincides with the measure $\mu$.
In other words, for every open set $A\subset\R^d$ there holds
\bel{chi1}\mu(A)~=~\hbox{meas}\Big( \{ \theta\in \Theta\,;~~\chi(\theta,T(\theta))\in A\bigr\}\Big).\eeq
\endi

Next, to define the corresponding transportation cost, one must 
take into account the fact that, if many paths go through the same pipe, their cost decreases.   With this in mind, 
given a point $x\in \R^d$ we first compute 
how many paths go through the point $x$.  
This is described by
\bel{chi}|x|_\chi~=~\meas\Big(\bigl\{\theta\in \Theta\,;~~\chi(\theta,t)= x~~~\hbox{for some}~~t\geq 0\bigr\}\Big).\eeq
We think of $|x|_\chi$ as the {\it total flux going through the
point $x$}.
\v
\begin{definition}
{\bf (irrigation cost).} 
For a given $\alpha\in [0,1]$,
the total cost of the irrigation plan $\chi$ is
\bel{TCg}
\E^\alpha(\chi)~\doteq~\int_\Theta\left(\int_0^{T(\theta)} \bigl|\chi(\theta,t)
\bigr|_\chi^{\alpha-1} \, dt\right)
d\theta.\eeq
The  {\bf $\alpha$-irrigation cost} of a measure $\mu$
is defined as 
\bel{Idef}\I^\alpha(\mu)~\doteq~\inf_\chi \E^\alpha(\chi),\eeq
where the infimum is taken over all admissible irrigation plans.
\end{definition}
\v
\begin{remark} {\rm In the case $\alpha=1$, the expression (\ref{TCg}) reduces to
$$
\E^\alpha(\chi)~=~\int_\Theta\left(\int_{\R_+} |\dot \chi_t(\theta,t)|\, dt\right)
d\theta~=~\int_\Theta[\hbox{total length of the path} ~\chi(\theta,\cdot)]\, d\theta\,.$$
Of course, this length is minimal if every path $\chi(\cdot,\theta)$
is a straight line, joining the origin with $\chi(\theta, T(\theta))$.  Hence
$$\I^\alpha(\mu)~=~\inf_\chi \E^\alpha(\chi)~=~\int_\Theta |\chi(\theta, T(\theta))|\, d\theta~=~\int |x|\, d\mu\,.$$

On the other hand, when $\alpha<1$, moving along a path which is traveled by few other particles
comes at a  high cost. Indeed, in this case the factor $\bigl|\chi(\theta,t)
\bigr|_\chi^{\alpha-1}$ becomes  large.   To reduce the total cost,  is thus convenient
that particles travel along the same path as far as possible.}
\end{remark}
\v
For the basic theory of ramified transport we refer 
to \cite{BS, MMS, X03}, or to the monograph \cite{BCM}. 
The following lemma provides a useful lower bound
to the transportation cost.  In particular, we recall that optimal irrigation plans
satisfy
\v
{\bf Single Path Property:}  If $\chi(\theta, \tau)=\chi(\theta',\tau')$ for some 
$\theta, \theta'\in\Theta$ and 
$0<\tau<\tau'$, then 
\bel{SPP}\chi(\theta, t)~=~\chi(\theta', t)\qquad\forall t\in [0, \tau].\eeq

\v
\begin{lemma}
For any positive Radon measure $\mu$ on  $\R^d$ and any $\alpha\in [0,1]$, one has
 \bel{Ilb}
\I^\alpha(\mu)~\geq~\int_0^{+\infty}\Big(\mu\bigl(\{x\in\R^d;~|x|\geq r\}
\bigr)\Big)^\alpha dr\,.\eeq
In particular, for every $r>0$ one has
 \bel{Irb}
\I^\alpha(\mu)~\geq~r\cdot \mu\bigl(\{x\in\R^d;~|x|\geq r\}
\bigr)^\alpha\,.\eeq
\end{lemma}
\v
{\bf Proof.}
Let $\chi: \Theta\times\R_+\mapsto \R^d$ be an optimal transportation
 plan for $\I^\alpha(\mu)$. 
For any given $t> 0$, let
$$\Theta_t~\doteq~ \Big\{\theta\in \Theta\,;\quad T(\theta)\geq t \Big\}$$
be the set of particles whose path has length $\ge t$.
By the Single Path Property (see Chapter~7 in \cite{BCM}),  
if 
$$\chi(\theta, \tau)~ = ~\chi(\tilde \theta, \tilde \tau),$$
for some $\theta,\tilde\theta\in \Theta$ and $0\leq \tau\leq\tilde\tau$, 
then 
\bel{spp}\chi(\theta, t) ~=~\chi(\tilde\theta,t)\qquad\forall t\in [0,\tau].\eeq
As a consequence, if $t\leq T(\theta)$, then
\bel{c<} \bigl|\chi(\theta, t)\bigr|_{\chi}~\leq~\meas(\Theta_t).\eeq

In addition, since all particles travel with unit speed, we have the obvious implication
$$
x~=~\chi(\theta, T(\theta))\qquad\implies\qquad T(\theta)~\geq~|x|,$$
hence
\bel{tt}
\meas(\Theta_t)~\geq~\mu\Big(\bigl\{ x\,;~~|x|\geq t\bigr\}\Big).\eeq
Always relying on the optimality of  $\chi$,  by (\ref{c<}) and (\ref{tt}) we conclude
$$\bega{rl} 
\I^\alpha(\mu)&=~\ds
\E^\alpha(\chi)~\doteq~\int_\Theta\left(\int_{\R_+} \bigl|\chi(\theta,t)
\bigr|_\chi^{\alpha-1} \cdot |\dot \chi(\theta,t)|\, dt\right)
d\theta\\[4mm]
&\ds =~\int_\Theta\left(\int_0^{T(\theta)} \bigl|\chi(\theta,t)
\bigr|_\chi^{\alpha-1}\, dt\right)d\theta ~\ge~\int_\Theta\left(\int_0^{T(\theta)}
\bigl[\meas(\Theta_t)\bigr]^{\alpha-1} \,dt\right)d\theta\\[4mm]
&\ds =~\int_0^{+\infty} \left(\int_{\{T(\theta)\geq t\}} \bigl[\meas(\Theta_t)\bigr]^{\alpha-1}
\,d\theta\right)dt~=~\int_0^{+\infty} \bigl[\meas(\Theta_t)\bigr]^\alpha
\,dt\\[4mm]
&\geq~\ds \int_0^{+\infty} \Big[\mu\bigl(\{x\in\R^d;~|x|\geq t\}
\bigr)\Big]^\alpha\,dt\,.
\enda $$
This proves (\ref{Ilb}).   The inequality (\ref{Irb}) follows immediately.
\endproof
\v
\section{Existence of optimal branch configurations, without constraint on the total mass}
\label{s:2}
\setcounter{equation}{0}
In this section we study a problem related to  the optimal shape of tree branches.
\begi
\item[{\bf (OPB)}] {\bf Optimization Problem for  Branches.} {\it  Given a function 
$\eta\in \L^1(S^{d-1})$ and constants $\alpha\in [0,1], c>0$,
\bel{OB}\hbox{maximize:} ~~S^\eta(\mu)-c\I^\alpha(\mu),\eeq
among all positive Radon measures  $\mu$, supported on closed the half space
\bel{barR}
\ov\R^d_+~\doteq~\bigl\{(x_1,\ldots, x_d)\,;~~x_d\geq 0\bigr\},\eeq
 without any constraint on the total mass.}
\endi

In \cite{BSu} the existence of an optimal solution to the problem
(\ref{OB})
was proved under a constraint on the total mass of the measure $\mu$, namely
$$
\mu(\R^d)~\leq~\kappa_0.$$
Our present goal is to prove  the existence of an optimal  solution of (\ref{OB}) 
without any  constraint.

Throughout the following, it will be natural to assume 
\bel{aas}
1-{1\over d-1}~\doteq~\alpha^*~<~\alpha~\leq~1.\eeq
Indeed, if a measure $\mu$ is supported on a set whose $(d-1)$-dimensional measure is zero,
then $\S^\eta(\mu)=0$.  On the other hand, if $\alpha<\alpha^*$, then any set 
with positive $(d-1)$-dimensional measure cannot be irrigated.  
Therefore, for $\alpha<\alpha^*$ the optimization problem (\ref{OB}) becomes trivial: the
zero measure is already an 
optimal solution.  We can now state 
the main result of this section.

\begin{theorem}\label{t:21}
 Suppose  that
$d\geq 3$ and $\alpha>\alpha^*$, as in (\ref{aas}). Then the 
unconstrained optimization problem (\ref{OB}) admits an optimal solution $\mu$, with bounded 
support and bounded total mass.
\end{theorem}

{\bf Proof.}   Following the direct method in the Calculus of Variations,  we consider a maximizing sequence of measures
$(\mu_k)_{k\geq 1}$.  While each $\mu_k$ is a bounded positive measure, at this stage
we cannot exclude the possibility that $\mu_k(\R^d)\to +\infty$.   By showing that all 
measures $\mu_k$ are uniformly bounded and have uniformly bounded support, we
shall be able to select a subsequence, weakly converging to an optimal solution.
The proof is given in several steps.
\v
 {\bf 1.} 
As a first step, we claim that the irrigation costs $\I^\alpha(\mu_k)$
are uniformly bounded.

Indeed,  given a radius $r>0$, we can decompose any measure $\mu$  as a sum
\bel{md1}\mu~=~\mu^-_r + \mu^+_r
~\doteq~\chi_{\{x\leq r\}}\cdot \mu + \chi_{\{x> r\}}\cdot \mu
\,.
\eeq
Here $\chi_A$ denotes the characteristic function of a set $A\subset\R^d$.
Calling $\omega_{d-1}$ the volume of the unit ball in $\R^{d-1}$,  using 
(\ref{Sbo1})-(\ref{Sbo2})  and then (\ref{Irb}),  the sunlight functional can now be bounded as
\bel{sf1} \bega{rl}\S^\eta(\mu)&\leq~\S^\eta(\mu_r^-)+\S^\eta(\mu_r^+)\\[3mm]
&\leq~\|\eta\|_{\L^1}\cdot \omega_{d-1} r^{d-1} + \|\eta\|_{\L^1}\cdot \mu\Big(\{x\,;~|x|>r\}\Big)
\\[3mm]
&\le~\ds \|\eta\|_{\L^1}\cdot \left[\omega_{d-1} r^{d-1} + \left( {I^\alpha(\mu)\over r}
\right)^{1/\alpha}\right].
\enda\eeq
In the above inequality, the  radius $r>0$ is arbitrary. 
In particular, we  can choose $r$ such that 
$$\omega_{d-1}r^{d-1}~=~\Big(\frac{\I^\alpha(\mu)}{r}\Big)^{1/\alpha}.$$
This choice yields
\bel{ria}\omega_{d-1}^\alpha r^{1+\alpha(d-1) }~=~I^\alpha(\mu),
\qquad\qquad 
r~=~\left( I^\alpha(\mu)\over \omega_{d-1}^\alpha\right)^{1\over 1+\alpha(d-1)}.\eeq
Inserting (\ref{ria}) in (\ref{sf1}) one obtains the a priori bound
\bel{Seb}
\S^\eta(\mu)~\leq~C_0\, \Big(\I^\alpha(\mu)\Big)^{d-1\over 1+\alpha(d-1)},
\eeq
for some constant $C_0$ depending only on $\alpha, d$, and $\|\eta\|_{\L^1}$.

In connection with the original problem (\ref{OB}), this implies
\bel{Ibb}
\S^\eta(\mu)-c\I^\alpha(\mu)~\leq~C_0\, \Big(\I^\alpha(\mu)\Big)^{d-1\over 1+\alpha(d-1)}
-c\I^\alpha(\mu).\eeq
We now observe that 
the assumption (\ref{aas}) is equivalent to
$${d-1\over 1+\alpha(d-1)}~<~1.$$
Therefore, by (\ref{Ibb}) there exists a constant $\kappa_1$ large enough so that
\bel{k*} \I^\alpha(\mu)~\geq~\kappa_1\qquad\implies\qquad \S^\eta(\mu)-c\I^\alpha(\mu)~\leq~0.\eeq
In the remainder of the proof, without loss of generality we shall seek a global 
maximum for the functional in (\ref{OB}) under the additional constraint
\bel{IBB}
\I^\alpha(\mu)~\leq~\kappa_1\,.\eeq
In turn, by (\ref{Seb}) one has a uniform bound 
\bel{SEB}
\S^\eta(\mu)~\leq~\kappa_2\eeq
for all $\mu$ satisfying (\ref{IBB}).
\v
{\bf 2.} Let $(\mu_k)_{k\geq 1}$ be a  maximizing sequence. 
In this step we construct a second   maximizing sequence $(\tilde \mu_k)_{k\geq 1}$ 
such that all measures $\tilde \mu_k$ are supported inside a fixed ball
$B_\rho$ centered at the origin with radius $\rho$.

Toward this goal, let  $\chi$ be an optimal irrigation plan for a measure $\mu$, as in (\ref{iplan}).
By (\ref{Irb}) and (\ref{IBB}), for any radius $r>0$  one has
\bel{mrb}
\mu\Big(\{x\in\R^d\,;~~|x|\geq r\}\Big)~\leq~\left(\I^\alpha(\mu)\over r\right)^{1/\alpha}
~\leq~\left(\kappa_1\over r\right)^{1/\alpha}.\eeq

Consider two radii $0<r_1< r_2$. As in (\ref{md1}), we can decompose 
the measure $\mu$ as a sum:
\bel{md11}\mu~=~\mu^\flat + \mu^\sharp~\doteq
~\chi_{\{x\leq r_2\}}\cdot \mu + \chi_{\{x> r_2\}}\cdot \mu
\,.
\eeq
By possibly relabeling the set $\Theta=\Theta^\flat\cup\Theta^\sharp$, 
we can assume that 
\begi
\item $\chi^\flat:\Theta^\flat\times\R_+\mapsto\R^d$ is an irrigation plan for the measure 
$\mu^\flat$
\item $\chi^\sharp:\Theta^\sharp\times\R_+\mapsto\R^d$ is an irrigation plan for the measure 
$\mu^\sharp$.
\endi
Note that $\chi^\flat$ and $\chi^\sharp$ are not necessarily optimal.
If $\mu^\sharp$ is removed, by (\ref{mrb}) the difference 
in the gathered sunlight is 
\bel{uvw1}
\S^\eta(\mu^\flat+\mu^\sharp)-\S^\eta(\mu^\flat)~\leq~\S^\eta(\mu^\sharp)~\leq~
\|\eta\|_{\L^1} \cdot \mu^\sharp(\R^d)~\leq~\|\eta\|_{\L^1} \cdot 
\left(\kappa_1\over r_2\right)^{1/\alpha}.\eeq

On the other hand, 
by the Single Path Property (\ref{SPP}), for any  $x\in \R^d$ 
with $|x|\geq r_1$ one has
$$I^\alpha(\mu)~\geq~|x|_{\chi}^\alpha\cdot r_1\,.$$
Therefore
\bel{uvw0}
 |x|_{\chi}~\leq~ \left(\frac{\I^\alpha(\mu)}{r_1}\right)^{1/\alpha}~\leq~ \left(\frac{\kappa_1}{r_1}\right)^{1/\alpha}.
\eeq
We now  estimate the difference of the irrigation cost, if part of the measure is removed.
Two cases will be considered.

CASE 1: $0<\alpha <1$.
By (\ref{uvw0})  we can then choose $r_1$ large enough so that 
\bel{uvw11}
|x|~\geq~r_1\qquad\implies\qquad \alpha\,  |x|_\chi^{\alpha-1}~\geq~ 1.\eeq
 According to Proposition~4.8 in \cite{BCM}, the  cost of an irrigation plan $\chi$
 can be equivalently described as
\bel{equive}\E^\alpha(\chi)~=~\int_{\R^d}|x|_{\chi}^\alpha\, d\H^1(x).\eeq
where $\H^1$ denotes the 1-dimensional Hausdorff measure.
If $\chi$ is an optimal irrigation plan  for $\mu=\mu^\flat+\mu^\sharp$, then
\bel{oenway}
\bega{rl}
\I^\alpha(\mu^\flat+\mu^\sharp)&\ds=~\int_{\R^d}|x|_{\chi}^\alpha\, d\H^1(x)
~=~\int_{\R^d}\Big(|x|_{\chi^\sharp}+|x|_{\chi^\flat}\Big)^\alpha\, d\H^1(x)\\[4mm]
 &\geq\ds~\int_{\R^d}\Big(|x|_{\chi^\flat}^\alpha+\alpha\bigl( |x|_{\chi^\flat}+|x|_{\chi^\sharp}\bigr)^{\alpha-1}|x|_{\chi^\sharp}\Big)\, d\H^1(x)\\[4mm]
 &\geq\ds~\int_{\R^d}|x|_{\chi^{\flat}}^\alpha\, d\H^1(x)+\int_{\{ |x|\geq r_1\}}
 \alpha |x|_\chi^{\alpha-1}\cdot |x|_{\chi^\sharp}\, d\H^1(x)\\[4mm]
 &\geq\ds ~\I^\alpha(\mu^\flat)+\int_{\{ |x|\geq r_1\}}
 |x|_{\chi^\sharp}\, d\H^1(x)\\[4mm]
 &\geq~\I^\alpha(\mu^\flat)+(r_2-r_1)\,\mu^\sharp(\R^d).
\enda\eeq


We now choose $r_2$ large enough so that $c(r_2-r_1)\geq \|\eta\|_{\L^1}$.
 By the second inequality in (\ref{uvw1}) and (\ref{oenway}) it follows
 \bel{dif5}\S^\eta(\mu^\flat+\mu^\sharp)-\S^\eta(\mu^\flat)~\leq~c\Big(\I^\alpha(\mu^\flat+\mu^\sharp)-\I^\alpha(\mu^\flat)\Big).\eeq
 Let now  $(\mu_k)_{k\geq 1}$ be a maximizing sequence.   We decompose each 
 measure as
\bel{mdk}\mu_k~=~\mu^\flat_k + \mu^\sharp_k~\doteq
~\chi_{\{x\leq r_2\}}\cdot \mu_k + \chi_{\{x> r_2\}}\cdot \mu_k
\,.
\eeq
By (\ref{dif5}), the sequence $(\mu_k^\flat)_{k\geq 1}$ is still a maximizing sequence, where all measures are supported inside the fixed ball $B_{r_2}\,$. 

CASE 2: $\alpha = 1$.   In this case we simply choose 
\bel{r2}r_2~=~{1\over c}\,\|\eta\|_{\L^1}\,.\eeq
 In connection with the decomposition (\ref{mdk}), we have
$$\bega{l}\S^\eta(\mu^\flat+\mu^\sharp)-\S^\eta(\mu^\flat)~\leq~
\S^\eta(\mu^\sharp)~\leq~\|\eta\|_{\L^1}\cdot \mu^\sharp(\R^d)\\[4mm]\qquad =~
r_2 c \,\mu^\sharp(\R^d)~\leq~c \I^1(\mu^\sharp)~=~c\Big(\I^1(\mu^\sharp)- \I^1(\mu^\flat+\mu^\sharp)\Big).\enda$$
 Again, this shows that  $(\mu_k^\flat)_{k\geq 1}$ is a maximizing sequence, where all measures are supported inside the ball $B_{r_2}\,$. 
\v
{\bf 3.} In this step, relying on the assumption that the space dimension is $d\geq 3$,
 we prove the existence of  a maximizing sequence $(\tilde \mu_k)_{k\geq 1}$ 
 with uniformly bounded total mass.

Indeed,  let $\mu$ be any  measure with $\I^\alpha(\mu)\leq \kappa_1$.
For any integer $j$, 
consider the radius  $r_j\doteq 2^{-j}$ and the spherical shell
\bel{Vj}V_j~\doteq ~\bigl\{ x\in\R^d\,;~~r_{j+1}<|x|\leq r_j\bigr\}.\eeq
Moreover, call 
$$\mu_j   ~\doteq~\chi_{\strut V_j}\cdot\mu$$
the restriction of the measure $\mu$ to the set $V_j$.  For every $j\geq 1$ we then have
\bel{muj} \S^\eta(\mu)-
\S^\eta(\mu-\mu_j)~\leq~\S^\eta(\mu_j)~\leq~\|\eta\|_{\L^1}\cdot \omega_{d-1} \,r_j^{d-1}.
	\eeq 
We now estimate the difference in the irrigation costs.
By (\ref{uvw0}), for every $x\in\R^d$ one has
\bel{mm} |x|_{\chi}~\leq~\left({\kappa_1\over |x|}\right)^{1/\alpha},\eeq
hence
$$\min\left\{|z|_\chi^{\alpha-1}\,;~~|z|\geq r_{j+2} \right\}~\geq~
\kappa_1^{\alpha-1\over\alpha}
\cdot r_{j+2}^{{1\over\alpha}-1}~\doteq~\kappa_3 \,r_j^{{1\over\alpha}-1}, $$
for a suitable constant $\kappa_3$.
This implies
\bel{Idif}
\I^\alpha(\mu)-\I^\alpha(\mu-\mu_j)
~\geq~{r_{j}\over 4}\cdot \int_{V_j} \kappa_3\, r_j^{{1\over\alpha}-1} d\mu
~=~ {\kappa_3\over 4}\, r_j^{1/\alpha}\cdot\mu(V_j).\eeq
Comparing (\ref{muj}) with (\ref{Idif})
 we see that, if 
\bel{jbig} {\kappa_3\over 4}\, r_j^{1/\alpha}\cdot\mu(V_j)~\geq~\|\eta\|_{\L^1}\cdot \omega_{d-1} r_j^{d-1}, \eeq
then the difference $\S^\eta(\mu)-c\I^\alpha(\mu)$ will increase
if we  remove from $\mu$ all the mass located inside $V_j$.

We can repeat the above procedure, removing from $\mu$ the mass contained in all 
regions $V_j$ such that (\ref{jbig}) holds.  More precisely,
let $J$ be the set of all integers $j\geq 0$ for which (\ref{jbig}) holds, and consider the measure 
\bel{tmdef}
\tilde\mu~\doteq~\mu - \sum_{j\in J}\mu_j\,.\eeq
 By the previous analysis, 
 \bel{tmb}
  \S^\eta(\tilde \mu)-c\I^\alpha(\tilde\mu)~\geq~\S^\eta(\mu)-c\I^\alpha(\mu).\eeq
Moreover, by (\ref{jbig}) we have the implication
$$j\notin J\qquad\implies\qquad \mu(V_j)~\leq~\|\eta\|_{\L^1}\cdot {4\omega_{d-1}\over \kappa_3} \, r_j^{d-1-{1\over\alpha}}~\doteq~\kappa_4 \, r_j^{d-1-{1\over\alpha}}.$$
The total mass of $\tilde\mu$ can thus be estimated by
\bel{tmm}\bega{rl}
\tilde\mu(\R^d)&\ds=~\mu\Big(\{x\,;~~|x|>1\}\Big)+\sum_{j\notin J} \mu(V_j)
~\leq~\mu\Big(\R^d\setminus B_1\Big) + \sum_{j\geq 0}\kappa_4 \, r_j^{d-1-{1\over\alpha}}
\\[4mm]
&\ds =~\mu\Big(\R^d\setminus B_1\Big)+ \kappa_4 \,
\sum_{j\geq 0} 2^{-j\Big(d-1-{1\over\alpha}\Big)}~<~+\infty,
\enda\eeq
 provided that $d-1-{1\over\alpha}>0$.
This is indeed the case if $\alpha$ satisfies (\ref{aas}) and $d\geq 3$.
	\v	
{\bf 4.} By the previous steps, we can choose a maximizing sequence 
$(\mu_k)_{k\geq 1}$ such that all measures $\mu_k$ have  uniformly 
bounded total mass and are supported on a fixed ball. By possibly taking a subsequence, 
we achieve the weak convergence $\mu_k\wto\mu$ for some
bounded measure $\mu$. By the upper semicontinuity of sunlight functional $\S^\eta$ 
proved in \cite{BS} and by the lower semicontinuity of the irrigation cost $\I^\alpha$, 
see \cite{BCM, MMS}, this limit measure $\mu$ provides a solution 
to the optimization problem (\ref{OB}).
\endproof

\subsection{The case $d=2$.}

In dimension $d=2$ we have  $d-1-{1\over\alpha}\leq 0$ for all $\alpha\leq 1$,
hence the estimate (\ref{tmm}) on the total mass breaks down. 
We develop here a different approach,
which is valid for
\bel{abig}
{\sqrt 5 -1\over 2}~<~\alpha~\leq~1\,.\eeq

\begin{figure}[ht]
\centerline{\hbox{\includegraphics[width=9cm]{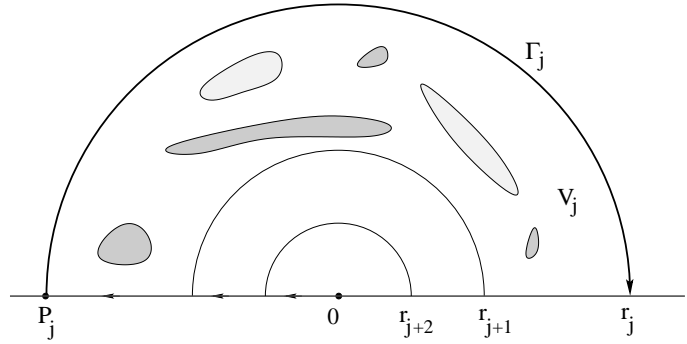}}}
\caption{\small  In dimension $d\geq 3$,  if  $\mu(V_j)$ is large, then we can
increase the 
payoff (\ref{OB}) by simply removing all the mass contained in the spherical shell $V_j$.
This idea is used in step {\bf 3} of  the proof of Theorem~\ref{t:21}.
In dimension $d=2$, if $\mu(V_j)$ is large, to  increase the 
payoff (\ref{OB}) we replace the measure 
$\mu_j = \chi_{V_j}\cdot\mu$
with a new measure $\tilde \mu_j$
uniformly distributed over the half circumference $\Gamma_j$. Notice that $\tilde\mu_j$ can be 
irrigated by moving the water particles from the origin to $P_j$, and then along $\Gamma_j\,$.}
\label{f:ir16}
\end{figure}
\v
\begin{theorem}\label{t:22}
If $d=2$  and $\alpha$ satisfies (\ref{abig}), then the 
unconstrained optimization problem (\ref{OB}) admits an optimal solution $\mu$,  with bounded 
support and bounded total mass.
\end{theorem}

Indeed, repeating the steps {\bf 1 - 2} in the proof of the Theorem~\ref{t:21}, we obtain
a maximizing sequence $(\mu_k)_{k\geq 1}$ of positive measures with uniformly bounded
support. Moreover, the 
irrigation costs $\I^\alpha(\mu_k)$ remain uniformly bounded.

In order to achieve a uniform bound on the total mass $\mu_k(\R^d)$,
an auxiliary result is needed.

\begin{lemma}\label{l:3.5} Let $\alpha$ satisfy (\ref{abig}) and let $\kappa_1>0$ be given.
Then  there exists an integer $j^*$ and an exponent $\ve>0$ 
such that the following holds.  Given any bounded measure $\mu$ with $\I^\alpha(\mu)\leq \kappa_1$, there exists a second measure $\tilde\mu$ satisfying (\ref{tmb}) and such that, setting $r_j\doteq 2^{-j}$,
\bel{mje}\tilde\mu\Big(\bigl\{ x\in\R^2\,;~~r_{j+1}<|x|\leq r_j\bigr\}\Big)
~\leq~2^{-\ve j}
\qquad\qquad\forall j\geq j^*.\eeq
\end{lemma}
\v
{\bf Proof.} {\bf 1.} If (\ref{abig}) holds, we can find $0<\ve<\beta<1$ such that 
\bel{ap1}\alpha\beta +1~>~\ve + {1\over\alpha}\,.\eeq
Let $\mu_j$ be the restriction of the measure $\mu$ to the spherical shell $V_j$ 
defined at  (\ref{Vj}).
 Moreover, let $\tilde\mu_j$ be the positive measure with total mass
 $$\tilde\mu_j(\R^d)~=~\pi r_j^{\,\beta},$$
 uniformly distributed on the half circumference
 $$\Gamma_j~\doteq~\{ x=(x_1,x_2)\in\R^2\,;~~|x|= r_j\,,~x_2>0\}.
$$
As shown in Fig.~\ref{f:ir16}, 
there is a simple irrigation plan $\chi$ for $\tilde\mu_j$.   Namely, we can 
first move all water particles on a straight line from the origin to the point $P_j = (-r_j, 0)$,
then from $P_j$ to all points along the half circumference $\Gamma_j$.
The total cost of this irrigation plan satisfies
\bel{Imj}\bega{rl}
\E^\alpha(\chi)&
\leq~[\hbox{total mass}]^\alpha\times [\hbox{maximum length traveled}]\\[3mm]
&\leq~(\pi r_j^{\beta})^\alpha\cdot (\pi+1) r_j.\enda\eeq
Therefore, 
the minimum irrigation cost for $\tilde\mu_j$ satisfies
\bel{Imj2}
\I^\alpha(\tilde\mu_j)\leq~2 \pi^{\alpha+1}\, r_j^{\alpha\beta+1}.\eeq
On the other hand, assuming $\mu(V_j)\geq r_j^{\,\ve}$, by (\ref{Idif}) we have
\bel{ID2}
\I^\alpha(\mu)-\I^\alpha(\mu-\mu_j)~\geq~{\kappa_3\over 4} r_j^{\ve+{1\over\alpha}}.\eeq
By (\ref{ap1}), for all $r_j$ small enough it follows
\bel{CI3}\Big[ \I^\alpha(\mu)-\I^\alpha(\mu-\mu_j)\Big] -
\I^\alpha(\tilde\mu_j) ~\geq~{\kappa_3\over 8} r_j^{\ve + {1\over\alpha}}.\eeq
\v

\begin{figure}[ht]
\centerline{\hbox{\includegraphics[width=9cm]{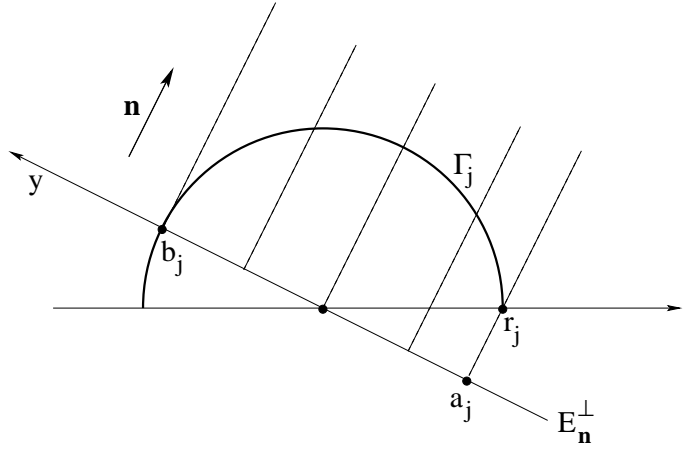}}}
\caption{\small 
Let $\tilde \mu_j$ be the measure supported on the half circumference $\Gamma_j$, 
with constant density $r_j^{\beta-1}$ w.r.t.~1-dimensional measure. 
Then, for any unit vector $\bfn$, the projected measure $\pi_\bfn\tilde\mu_j$ has density 
$\geq~r_j^{\beta-1}$ on $[a_j, b_j]=\pi_\bfn(\Gamma_j)$.}
\label{f:ir17}
\end{figure}

{\bf 2.} Next, we estimate how  the sunlight functional  changes
if we replace $\mu_j$ by $\tilde \mu_j$.
We claim that 
\bel{SD2}\bega{l}
\S^\eta(\mu)-S^\eta(\mu - \mu_j+\tilde\mu_j)\\[3mm]
\qquad \le~[\hbox{total amount of light hitting $V_j$}]
- [\hbox{light captured by $\tilde\mu_j$}]\\[3mm]
\qquad \le~\|\eta\|_{\L^1}\cdot \exp\bigl(-r_j^{\beta-1}\bigr).\enda\eeq
Indeed, consider any unit vector $\bfn\in S^1$. 
As shown in Fig.~\ref{f:ir17}, let $[a_j, b_j]=\pi_\bfn(\Gamma_j)$
 be the 
perpendicular projection of $\Gamma_j$ on the orthogonal subspace 
$E^\perp_\bfn$.
By construction, the projected measure $\pi_\bfn\tilde\mu_j$ is absolutely 
continuous w.r.t.~1-dimensional Lebesgue measure on $E^\perp_\bfn$. 
Its density $\Phi_j^\bfn$ satisfies
$$\left\{ \bega{rll}\Phi_j^\bfn(y)&\geq~r_j^{\beta-1}\quad &\hbox{if}~~y\in [a_j, b_j],\\[2mm]
\Phi_j^\bfn(y)&=~0\quad &\hbox{if}~~y\notin [a_j, b_j].\enda\right.$$
For $j\geq 1$ we thus have 
\bel{SD3}\bega{l}
[\hbox{total amount of light   hitting $V_j$ in the direction $\bfn$}]\\[3mm]
\qquad\qquad 
- [\hbox{light parallel to $\bfn$ captured by $\tilde\mu_j$}]\\[3mm]
\ds\qquad \le~(b_j-a_j) - \int_{a_j}^{b_j} \Big( 1- e^{-\Phi_j^\bfn(y)}\Big)\, dy
\\[3mm]
\ds\qquad 
\le~\int_{a_j}^{b_j} \exp\bigl(- r_j^{\beta-1}\bigr)\, dy~\leq~\exp\bigl(- r_j^{\beta-1}\bigr),
\enda\eeq
because $b_j-a_j\leq 1$.
\v
{\bf 3.} We now observe that,
since $\beta<1$, when $j\geq j^*$ is sufficiently large  the right hand side of 
(\ref{SD3}) is smaller than the right hand side of (\ref{CI3}).    By possibly 
choosing a larger $j^*$, we can also assume that
\bel{peb}\pi r_j^\beta ~<~r_j^\ve\qquad\forall j\geq j^*.\eeq
Defining the set of indices
$$J~\doteq~\bigl\{ j\geq j^*\,;~~\mu(V_j)~>~r_j^\ve\bigr\},$$
we claim that  the modified measure
\bel{tmde}\tilde\mu~\doteq~\mu + \sum_{j\in J} (\tilde\mu_j-\mu_j).\eeq
satisfies all conclusions of the lemma.
Indeed,  from (\ref{CI3}) and (\ref{SD2})-(\ref{SD3}) it follows that
$\tilde\mu$ achieves a better payoff than $\mu$, i.e.~(\ref{tmb}) holds. 
In addition, the bounds
(\ref{mje}) on the total mass follow from (\ref{peb}).  
\endproof
\v
We observe that  (\ref{mje}) implies an a priori bound on the total mass
$\tilde \mu\Big(\bigl\{ x\in\R^2\,;~~|x|\leq r_{j^*}\bigr\}\Big)$.
On the other hand,   a bound on
$\tilde \mu\Big(\bigl\{ x\in\R^2\,;~~|x|\geq r_{j^*}\bigr\}\Big)$
is already provided by (\ref{mrb}).

Thanks to the above lemma, the proof of Theorem~\ref{t:22} is now straightforward.
Indeed, by Lemma~\ref{l:3.5}  we can replace each $\mu_k$ 
by a new measure $\tilde\mu_k$ we can one can construct a maximizing sequence 
of measures with uniformly bounded support and uniformly bounded total mass. 
Taking a weak limit, the existence of an 
optimal solution can thus be proved using the upper semicontinuity 
of $\S^\eta$ and the lower semicontinuity of $\I^\alpha$, as in 
\cite{BSu}.
\v
\section{Optimal root configurations, without size constraint}
\label{s:25}
\setcounter{equation}{0}
In this section we study the optimal shape of tree roots.
\begi
\item[{\bf (OPR)}] {\bf Optimization Problem for Roots.} {\it
\bel{OR}\hbox{maximize:} ~~\H(u,\mu)-c\I^\alpha(\mu),\eeq
subject to
\bel{rm2}\left\{\bega{rl}
\Delta u + f(u) - u\mu&=~0,\qquad\qquad x\,\in\, \R^d_-\,\doteq\, 
\bigl\{(x_1,\ldots, x_d)\,;~~x_d<0\bigr\},\\[3mm]
u_{x_d}&=~0,\qquad\qquad  x_d\,=\,0.\enda\right.\eeq
Here $\mu$ is a positive measure  concentrated on the set 
\bel{om0}
\Omega_0~\doteq~\bigl\{(x_1,\ldots, x_d)\not= (0,\ldots, 0)\,;~~x_d\leq 0\bigr\},\eeq
 without any constraint on its total mass.}
\endi
We recall that
\bel{HF4}\H(u,\mu)~=~\int_{\{x_d\,\leq\, 0\}} u \,d\mu~=~\int_{\R^d_-} f(u)\, dx\eeq
is the harvest functional introduced at (\ref{HF1}),
while $\I^\alpha(\mu)$ is the minimum irrigation cost defined at (\ref{Idef}).

As in Section~2, we assume that the function $f$  satisfies all conditions  
in (\ref{fp1}).
In order to construct an optimal solution, we consider
a maximizing sequence $(u_k, \mu_k)_{k\geq 1}$.   By suitably 
adapting the arguments used 
in the previous section, we will prove
a priori bounds on the total irrigation costs $\I^\alpha(\mu_k)$ 
and on the total harvesting payoffs $\H(u_k,\mu_k)$.
Our first lemma shows that the total harvest achieved by a measure supported on a closed
ball $\ov B_{\rho}$, centered at the origin with a large radius $\rho$, grows at most like $\rho^d$.
\v
\begin{lemma}\label{l:3.1} Let $f$ satisfy the assumptions (A1).
Then there exists a constant $C_f$  such that the following holds.
For any $\rho\geq 1$, if $\mu$ is a positive measure supported inside the closed ball
$\ov B_{\rho}$, then for any solution $u$ of (\ref{rm2}) one has 
\bel{Hub}\H(u,\mu)~\leq~C_f \,\rho^d.\eeq
\end{lemma}

\begin{figure}[ht]
\centerline{\hbox{\includegraphics[width=5cm]{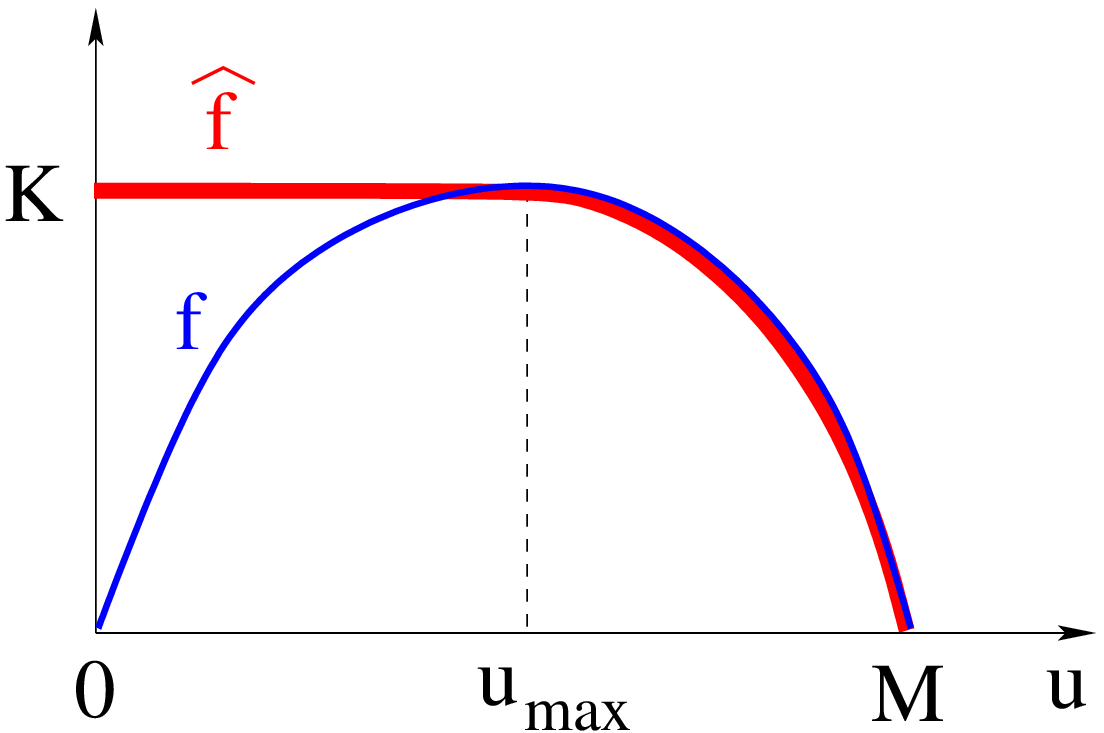}}\quad\hbox{\includegraphics[width=10cm]{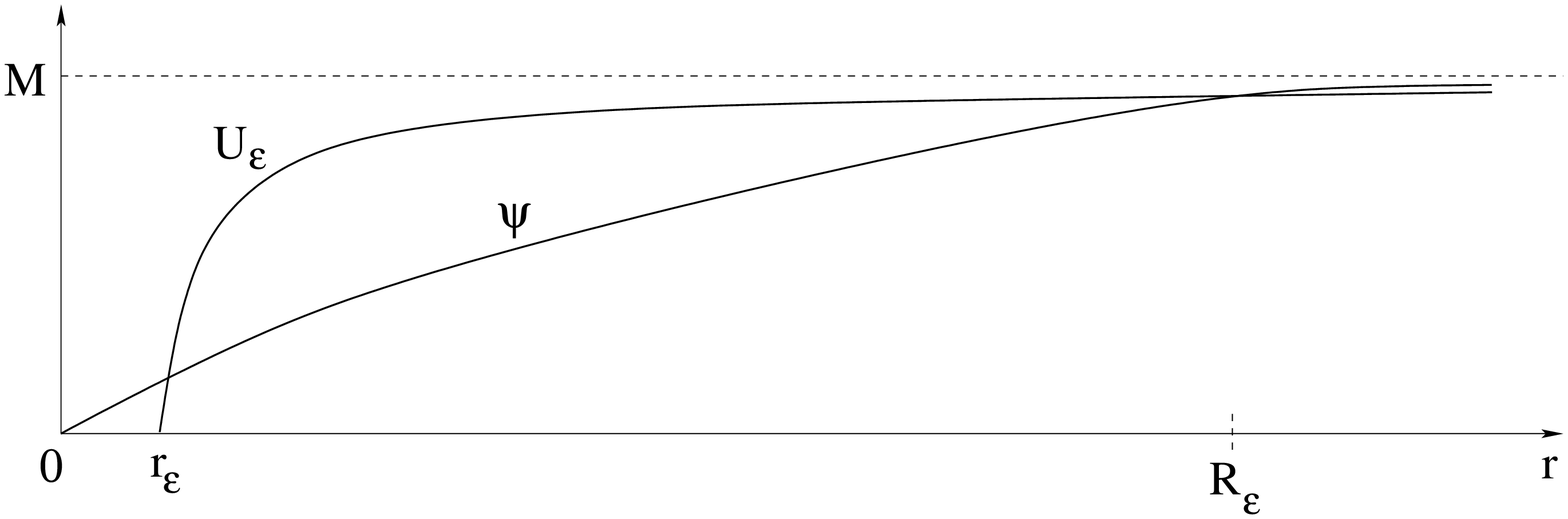}}}
\caption{\small Left: a function $f$ satisfying the assumptions
(\ref{fp1}) and the corresponding function $\Hat f$ in (\ref{hfdef}).
Right: a subsolution obtained by patching  together the functions
$U_\ve$ and $\psi$. }
\label{f:ir37}
\end{figure}

{\bf Proof.} {\bf 1.} As shown in Fig.~\ref{f:ir37}, right,
 let $\psi=\psi(r)$ be the solution to the ODE
\bel{subs}\psi''(r) 
+f(\psi(r))~=~0,\qquad\qquad r>0,\eeq
\bel{sub2}
\psi(0)~=~0,\qquad\qquad \lim_{r\to +\infty} \psi(r)~=~M.\eeq
We claim that $\psi$ is a monotonically increasing function such that 
$$\psi(r)~\to~M\qquad \hbox{as}\qquad r\to +\infty,$$
with an exponential rate of convergence.

Indeed, let $F(s)=\int_0^s f(\xi)\, d\xi$.   Then, for any solution of (\ref{subs}),
the energy
\bel{encst}E(r)~\doteq~{\bigl(\psi'(r)\bigr)^2\over 2 }+ F\bigl(\psi(r)\bigr)\eeq
is constant.
The second limit in (\ref{sub2}) implies that $E=F(M)$.
We thus obtain the ODE
\bel{newode}\psi'(r)~=~\sqrt{2F(M)-2F(\psi(r))}.\eeq
Since $f(M)=0$ and $f'(M)<0$,
for $\psi\in [0,M]$ one has
\bel{engi}F(M)-F(\psi)~\geq~\gamma\,(M-\psi)^2,\eeq
for some constant $\gamma>0$ depending only on $f$ itself. Therefore,
$$\psi'(r)~>~\sqrt {2 \gamma}(M-\psi)$$
 We thus conclude that the solution
$\psi$ of (\ref{subs})-(\ref{sub2}) satisfies
\bel{expc}\psi(r)~\geq~M(1 - e^{-\sqrt{2\gamma}\, r}).\eeq
Therefore
\bel{fpb}
f(\psi(r))~\leq~ f^{'}(M)(\psi(r)-M)~\leq~C_1 e^{-\sqrt{2\gamma}\, r}\, ,\qquad  C_1~\doteq~ -f^{'}(M)M>0\, .\eeq

\v
{\bf 2.} Let $u$ be a solution to (\ref{emu})-(\ref{NBC}), where the measure $\mu$ is supported
on the ball $\ov B_{\rho}$.  We claim that $u\geq v$, where $v$ is the radially symmetric function defined by
\bel{122}v(x)~=~\left\{\bega{cl}
\psi(|x|-\rho)\quad &\hbox{if}~~ |x|\geq \rho,\\[3mm]
0\quad &\hbox{if}~~ |x|< \rho.\enda\right.\eeq
Indeed, for $|x|>\rho$, by (\ref{122})  and (\ref{subs}) one has
	\bel{121}\bega{ll} \Delta v(x)+f(v)&\ds=~\psi''(|x|-\rho) +\frac{d-1}{|x|}\psi'(|x|-\rho)+
	f\bigl(\psi(|x|-\rho)\bigr)\\[3mm]
	&\ds =~\frac{d-1}{|x|}\psi'(|x|-\rho)~\geq~0,
	\enda\eeq
showing that $v$ is a lower solution on the region where $|x|>\rho$.
Hence $u(x)\geq v(x)$ for all $x\in\R^d$.
\v
{\bf 3.} Since $u\geq v$,
an upper bound on the total harvest is now provided by:
\bel{hub1}
\H(u,\mu)~=~\int_\Omega f(u(x))\, dx~\leq~{\omega_{d-1}\over 2}\,\int_0^{+\infty} r^{d-1}
\Hat f(\psi(r-\rho))\, dr,\eeq
where (see Fig.~\ref{f:ir37}, left)
\bel{hfdef}\Hat f(s)~=~\max\bigl\{f(\xi)\,;~~\xi\in [s,M]\bigr\}~=~\left\{ 
\bega{cl} f(s)\quad &\hbox{if} ~~s\geq u_{max}\,,\\[2mm]
K\quad &\hbox{if} ~~s\leq  u_{max}\,.\enda\right.\eeq
Here $u_{max}\in \,]0, M[\,$ is the unique point at which the function $f$ attains its maximum.
\v
{\bf 4.}	By the previous steps, the solution $\psi$ of (\ref{subs})-(\ref{sub2}) is a monotonically increasing function converging to $M$ as $s\to +\infty$. 
We can thus find a radius $r^*>1$ large enough so that 
\bel{r*} \psi(r)~\geq ~u_{max}\qquad\qquad\forall r\geq r^*.\eeq
Indeed, one can choose 
\bel{oldr}r^* ~\doteq~ -\frac{\ln(1 - \frac{u_{max}}{M})}{\sqrt{2\gamma}}\, .\eeq
Using (\ref{hfdef}),(\ref{fpb}) and performing  the variable change $r=\rho+s$, one finds
	\bel{123}\bega{rl}
\ds \int_0^\infty r^{d-1}\Hat f(\psi(r-\rho))dr&\ds =~\int_0^{\rho+r^*}r^{d-1}f(u_{max})\, dr
+\int_{\rho+r^*}^{+\infty} r^{d-1}f(\psi(r-\rho))\,dr\\[4mm]
 &\ds \leq~\frac{K}{d}(\rho+r^*)^d+\int_{\rho+r^*}^{+\infty} r^{d-1}C\, e^{-\sqrt{2\gamma}(r-\rho)}\, dr 
\\[4mm]
 &\ds \leq~\frac{K}{d}(\rho+r^*)^d+\rho^{d-1} 
 \int_{ r^*}^{+\infty} (1+s)^{d-1}C\, e^{-\sqrt{2\gamma}\, s}\, ds\\[4mm]
 &\leq~C_1 \rho^d+C_2\rho^{d-1}.
 \enda\eeq 
 Combining (\ref{123}) with (\ref{hub1}) one obtains the desired inequality (\ref{Hub}).
 \endproof
\v

The next lemma provides an estimate
on the total harvest achieved by a measure supported in a small
ball $\ov B_{\rho}$, as $\rho\to 0$.

\begin{lemma}\label{l:32} Let $f$ satisfy  (A1). Then there exists a positive
continuous function $\eta$, with  $\lim_{s\to 0+}\eta(s)=0$, such that the following holds.
Let $(u,\mu)$ be any solution of (\ref{rm2}).    If $\mu$ is supported on 
the closed  ball $\ov B_\rho$, then the total harvest satisfies
\bel{TH} 
\H(u,\mu)~\leq~\eta(\rho).\eeq
\end{lemma}

{\bf Proof.} {\bf 1.} 
Let $U=U(r)$ be a solution to
\bel{UE}
U''(r) +{d-1\over r} U'(r) + f(U)~=~0,\qquad\qquad \rho<r<+\infty,\eeq
\bel{UB}U(\rho)~=~0,\qquad\qquad \lim_{r\to +\infty} U(r)~=~M.\eeq
A lower bound on $U$ will be achieved by constructing a suitable subsolution.

Observing that $f(U)\geq 0$ and $U'>0$, such a subsolution can be obtained
by patching together a solution to
\bel{UE1} U''(r) +{d-1\over r} U'(r) ~=~0,\qquad\qquad \rho<r<r^*,\eeq
with a solution of 
\bel{UE2}
U''(r)  + f(U)~=~0,\qquad\qquad  r^*<r<+\infty.\eeq
As in the proof of the previous lemma, let  $\psi$ be the solution to
(\ref{subs}) and (\ref{sub2}).   In addition,
a solution of (\ref{UE1}) with boundary condition
\bel{U0}U(\rho)~=~0\eeq
is found in the form
\bel{SOL}
U(r)~=~\left\{ \bega{cl}\kappa\, (\ln r - \ln\rho)\qquad\qquad &\hbox{if}~~d=2,\\[3mm]
\kappa(\rho^{\,2-d} - r^{2-d})\qquad\qquad &\hbox{if}~~d\geq 3.\enda\right.\eeq
By linearity, here $\kappa$ can be any constant.
\v

{\bf 2.} 
To patch together the two solutions $U$ and $\psi$, we proceed as follows.
Recalling (\ref{hfdef}),  for any $\ve>0$, choose $R_\ve$ large enough so that 
\bel{ipe1}   R_\ve \,\bigl(\psi(R_\ve)-M\bigr)~<~\ve,\qquad\qquad \int_{R_\ve}^{+\infty} \Hat f(\psi(s))\, 
ds~\leq~\ve.\eeq
This is certainly possible because $f(M)=0$ and $\psi(s)\to M$ exponentially fast as $s\to +\infty$. 
\v
Next, we claim that there exists $r_\ve>0$ small enough and $\kappa_\ve>0$ so that
the function
\bel{Uedef}U_\ve(r)~\doteq~\left\{ \bega{cl}\kappa_\ve \, (\ln r - \ln r_\ve)\qquad\qquad &\hbox{if}~~d=2,\\[3mm]
\kappa_\ve(r_\ve^{\,2-d} - r^{2-d})\qquad\qquad &\hbox{if}~~d\geq 3,\enda\right.\eeq
satisfies (see Fig.~\ref{f:ir37}, right)
\bel{ipe6} U_\ve(R_\ve)~=~\psi(R_\ve),\qquad\qquad U'_\ve(R_\ve)~\leq~ \psi'(R_\ve),\eeq
\bel{ipe2} K r_\ve +\int_{r_\ve}^{R_\ve} \Hat f(U_\ve(s))\,ds~<~\ve.\eeq
Here $K$ is the maximum value of $f$, as in (\ref{fp1}). 

To prove our claim, having fixed $R_\ve$, for any $\rho>0$ we determine $\kappa_\ve$ so that
the function 
\bel{Uro}U_{\rho} (r)~\doteq~\left\{ \bega{cl}\kappa_\ve \, (\ln r - \ln \rho)\qquad\qquad &\hbox{if}~~d=2,\\[3mm]
\kappa_\ve(\rho^{\,2-d} - r^{2-d})\qquad\qquad &\hbox{if}~~d\geq 3,\enda\right.\eeq
satisfies   $U_\rho(R_\ve)=\psi(R_\ve)$.
As $\rho\to 0$, one now has
\bel{lu1}\lim_{\rho\to 0}
U'_\rho(R_\ve)~=~0,\qquad\qquad  \lim_{\rho\to 0} U_\rho(r)~=~\psi(R_\ve)\qquad\hbox{for any}~0<r\leq R_\ve.\eeq
This is proved by a direct computation.
When $d=2$ we have
$$\kappa_\ve~=~{\psi(R_\ve)\over \ln R_\ve -\ln\rho}\,,\qquad\quad
U_\rho(r)~=~\psi(R_\ve)\, {\ln r -\ln\rho\over \ln R_\ve -\ln\rho}\,,$$
$$U'_\rho(r)~=~{1\over r}\,{\psi(R_\ve)\over  \ln R_\ve -\ln\rho}\,.$$
Hence the limits in (\ref{lu1}) hold.
On the other hand, when $d\geq 3$ we have
$$\kappa_\ve~=~{\psi(R_\ve)\over R_\ve^{2-d} - \rho^{2-d}}\,,\qquad\quad
U_\rho(r)~=~\psi(R_\ve)\, {\rho^{2-d} - r^{2-d} \over \rho^{2-d} -R_\ve^{2-d} }\,,$$
$$U'_\rho(r)~=~{\psi(R_\ve)} { (d-2) r^{1-d} \over \rho^{2-d} -R_\ve^{2-d} }\,,$$
and the limits in (\ref{lu1}) again hold.
\v
Having determined $R_\ve$ according to (\ref{ipe1}), if we now choose
$r_\ve=\rho>0$ small enough, by  the first limit in 
(\ref{lu1}) it follows $U_\ve'(R_\ve)<\psi'(R_\ve)$.

Moreover, thanks to the second limit in (\ref{lu1}) and the first inequality in (\ref{ipe1}), 
by choosing $r_\ve>0$ small enough  we
also  achieve
$$Kr_\ve + \int_{r_\ve}^{R_\ve}  \Hat f(U_\ve(s))\, ds~<~\ve+R_\ve
 f(\psi(R_\ve))~\leq~\ve + R_\ve \,L_f (\psi(R_\ve)-M)~\leq~(L_f+1)\ve\,,$$
where $L_f$  denotes the Lipschitz constant of $f$.
Since $\ve>0$ was arbitrary, our claim is proved.
\v
{\bf 3.}
Let $u$ be a solution to (\ref{rm2}), where the measure $\mu$ is supported
on the closed ball $\ov B_{\rho}$.  By a comparison argument, we conclude  that $u\geq v$, 
where $v$ is the function defined by
$$v(x)~=~\left\{\bega{cl} 0\quad &\hbox{if}~~ |x|\leq r_\ve\,,\\[3mm]
U_\ve(r)\quad &\hbox{if}~~ r_\ve <|x|<R_\ve\,,\\[3mm]
\psi(r)\quad &\hbox{if}~~ R_\ve\leq |x|\,.
\enda\right.$$
By (\ref{ipe1})-(\ref{ipe2}),
an upper bound on the total harvest is now provided by
$$\bega{rl}
\H(u,\mu)&\ds =~\int_{\R^d_-} f(u(x))\, dx~\leq~{\omega_{d-1}\over 2}\,\int_0^{+\infty} r^{d-1}
\Hat f(v(r))\, dr\\[4mm]
&\ds 
=~{\omega_{d-1}\over 2} \left(\int_0^{r_\ve} \Hat f(0)\, dr +
\int_{r_\ve}^{R_\ve} U_\ve(r)\, dr +\int_{R_\ve}^{+\infty} \psi(r)\, dr\right) ~\leq ~\omega_{d-1}\, \ve\,.\enda$$
Since $\ve>0$ was arbitrary, this achieves the proof.
 \endproof
 \v
\begin{corollary}\label{c:33}
In the same setting as Lemma~\ref{l:32}, let $\mu$ be a positive measure on $\R^d$ 
and let $(u,\mu)$ be a solution to 
(\ref{rm2}). Then, for every $\rho>0$, one has
\bel{et1}
\int_{|x|\leq \rho} u\, d\mu~\leq~\eta(\rho),\eeq
where $\eta$ is the function in (\ref{TH}).
\end{corollary}

{\bf Proof.} Call $\mu^\rho\doteq\chi_{\strut \{|x|\leq \rho\}}\cdot\mu$ the restriction
of the measure $\mu$ to the closed ball of radius $\rho$.  Let $u^\rho\geq u$
be a corresponding solution of (\ref{rm2}).  Using  Lemma~\ref{l:32}
we now obtain
$$\int_{|x|\leq \rho} u\, d\mu~=~\int u\, d\mu^\rho~\leq~\int u^\rho\, d\mu^\rho~\leq~\eta(\rho).
$$
\endproof
  
Using Lemma~\ref{l:3.1}, we now prove that  an analogue of (\ref{k*}) holds also for the harvesting problem.
\v
\begin{lemma} Let $\alpha> 1-{1\over d}$. Under the assumptions (A1), 
there exists a constant $\kappa_1$ such that, for any positive bounded 
measure $\mu$ on $\ov\Omega\doteq \{ x= (x_1,\ldots, x_d)\in  \R^d\,;~~x_d\leq 0\}$,
\bel{H8} \I^\alpha(\mu)~\geq~\kappa_1\qquad\implies\qquad \H(u,\mu)- c\I^\alpha(\mu)~\leq ~0.\eeq
\end{lemma}

{\bf Proof.}  Given any radius $r\geq 1$, we can decompose the measure $\mu$ as
	\bel{rerfg}\mu~=~\mu_{r}^-+\mu_{r}^+~\doteq~\chi_{\{x\leq r\}}\cdot \mu + \chi_{\{x> r\}}\cdot \mu
	\,\eeq
	Let $u$ be a corresponding solution of  (\ref{rm2}).
 $0\leq u(x)\leq M$. 
Then there exists a solution $u^-$  to the same elliptic problem with $\mu$ replaced by 
$\mu_{r}^-$, such that 
\bel{uu-} 0~\leq~ u(x)~\leq~ u^-(x)\qquad\qquad\forall x\in\R^d_-\,.\eeq
Using (\ref{uu-}), and recalling that $\mu$ is concentrated on the domain 
$\Omega_0$ at (\ref{om0}), the harvest functional can be estimated by
\bel{wrpenn}
\bega{rl}\H(u,\mu)~&\ds \doteq~\int_{\Omega_0}u\, d\mu~=~\int_{\Omega_0}u\, d\mu_{ r}^-
+\int_{\Omega_0}u\, d\mu_{r}^+\\[4mm]
&\ds \leq~\int_{\Omega_0}u^-\, d\mu_{r}^-~+~M\mu_r^+(\Omega_0)\\[4mm]
 &\ds\leq~C_f\, r^d~+~M\Big(\frac{\I^\alpha(\mu)}{r}\Big)^{\frac{1}{\alpha}}.
 \enda\eeq
 Here the last inequality was obtained by applying  (\ref{Hub})  to the measure $\mu_r^-$
 and  
 (\ref{Irb}) to the measure $\mu^+_r$.

Next, assuming $\I^\alpha(\mu)$ sufficiently large, we can find  a radius $\rho\geq 1$ such that
\bel{stat}C_f\, \rho^{d}~=~M\Big(\frac{\I^\alpha(\mu)}{\rho}\Big)^{\frac{1}{\alpha}}.\eeq
Choosing $r=\rho$ in (\ref{wrpenn}), we obtain
\bel{sta1}\H(u,\mu)~\leq~C_0\Big(\I^\alpha(\mu)\Big)^{\frac{d}{1+\alpha d}}.\eeq
for some constant $C_0$ depending only on $\alpha,d,$ and $f$.

In connection with the original problem (\ref{OR}), this implies
\bel{aasse}\H(u,\mu)-c\I^\alpha(\mu)~\leq~C_0\Big(\I^\alpha(\mu)\Big)^{\frac{d}{1+\alpha d}}
-c\I^\alpha(\mu).\eeq

Assuming that $\alpha>1-\frac{1}{d}$, it follows that $\frac{d}{1+\alpha d}<1$.
Hence  by (\ref{aasse}) there exists a constant $\kappa_1$ large enough so that
(\ref{H8}) holds.
\endproof

\v
\begin{lemma}\label{l.35.1}
Let $\alpha> 1-{1\over d}$ and let  the assumptions (A1) hold.
Consider a maximizing sequence $(u_k, \mu_k)_{k\geq 1}$ for the functional 
in (\ref{OR}).   Then there exists another maximizing sequence $( \tilde u_k, \tilde\mu_k)_{k\geq 1}$
such that 
\bel{HIB}
\I^\alpha( \tilde \mu_k)~\leq~\kappa_1\,,\qquad\qquad 
\H( \tilde u_k, \tilde \mu_k)~\leq~\kappa_2\,,\eeq
for some constants $\kappa_1, \kappa_2$ and all $k\geq 1$.
Moreover, all measures $\tilde\mu_k$ are supported in a common ball
$\ov B_\rho$.
\end{lemma}

{\bf Proof.}  {\bf 1.} 
By (\ref{H8}), any maximizing sequence must satisfy
the first inequality in (\ref{HIB}). The second inequality then follows from
(\ref{sta1}).   
\v
{\bf 2.} 
It remains to prove the existence of a maximizing sequence with uniformly bounded support.
Toward this goal, let $\chi$ be an optimal irrigation plan for a measure $\mu$. By (\ref{Irb}) and (\ref{IBB}), for any radius $r>0$ one has
\bel{mrbaga}\mu\Big(\{x\in \R^d;~|x|\geq r \}\Big)~\leq~\Big(\frac{\I^\alpha(\mu)}{r}\Big)^{1/\alpha}~\leq~\Big(\frac{\kappa_1}{r}\Big)^{1/\alpha} .\eeq
Consider two radii $0<r_1<r_2$, where $r_1$ large enough such that (\ref{uvw11}) holds. As in (\ref{md111}), we can decompose the measure $\mu$ as a sum:
\bel{md111}\mu~=~\mu^{-} + \mu^+~\doteq~\chi_{\{x\leq r_2\}}\cdot \mu + \chi_{\{x> r_2\}}\cdot \mu\, .\eeq
By the same argument used in (\ref{oenway}), 
for $0<\alpha\leq 1$, one has
\bel{wyyy}\I^\alpha(\mu^- + \mu^+)~-~\I^\alpha(\mu^-)~\geq~(r_2 - r_1)\mu^+ (\Omega_0)\, ,\eeq
where $\Omega_0$ is the domain in (\ref{om0}).
\v
{\bf 3.}
We now estimate the decrease in the harvest functional, when  $\mu$ 
is replaced by $\mu^-$.   Let $u$ be a  solution of (\ref{rm2}), 
corresponding to the measure $\mu$.
Then there exists a solution $u^-$ to the same  problem, with $\mu$ replaced by $\mu^{-}$, such that 
\bel{vv-}0~\leq~u(x)~\leq~u^-(x)\qquad\quad \hbox{\forall} x\in \R^d_-\,.\eeq
Using (\ref{vv-}), the harvest functional can be estimated by
\bel{rpenst1}
\bega{rl}\H(u,\mu)~&\ds \doteq~\int_{\Omega_0}u\, d\mu~=~\int_{\Omega_0}u\, d\mu^-
+\int_{\Omega_0}u\, d\mu^+\\[4mm]
&\ds \leq~\int_{\Omega_0}u^-\, d\mu^-+\int_{\Omega_0}M\, d\mu^+\\[4mm]
&\ds =~\H(u^-,\mu^-)~+~M\mu^+(\Omega_0).
\enda\eeq
Taking $r_2$ large enough so that $c(r_2 - r_1) \geq M$, by (\ref{wyyy}) and (\ref{rpenst1}) it follows
\bel{qis}\H(u,\mu) - \H(u^-,\mu^-)~\leq~c\Big(\I^\alpha(\mu) - \I^\alpha(\mu^-)\Big).\eeq
\v
{\bf 4.} 
Let now $(u_k, \mu_k)_{k\geq 1}$ be a maximizing sequence. We decompose each measure as 
\bel{kmed}\mu_k ~=~\mu_k^- + \mu_k^+~\doteq~\chi_{\{x\leq r_2\}}\cdot \mu_k + \chi_{\{x> r_2\}}\cdot \mu_k\, .\eeq
Choose $u^-_k$ the corresponding solution to the same elliptic problem with $\mu_k$ replaced by $\mu_k^-$, such that 
\bel{vv+}0~\leq~u_k(x)~\leq~u_k^-(x)\qquad\quad \hbox{\forall } x\in \R^d_-\,.\eeq
By (\ref{qis}), $(u_k^-,\mu^-_k)_{k\geq 1}$ is still a maximizing sequence, where all measures are supported inside the closed ball $\ov B_{r_2}$.
\endproof
\v
The next lemma yields a more detailed estimate on the support of the optimal measure.

\begin{lemma}\label{l:36} Suppose $(u_k, \mu_k)_{k\geq 1}$
 is a maximizing sequence for the optimization problem (\ref{OR}), with irrigation costs
$\I^\alpha (\mu_k)\leq\kappa_1$ for all $k\geq 1$. Then there exists a second
maximizing sequence $(\tilde u_k, \tilde \mu_k)_{k\geq 1}$
such that
	\bel{lbdu} 
	\tilde\mu_k\Big(\bigl\{x\in \R^d\, ; \quad \tilde u_k(x)~<~C_0 |x|^{1/\alpha}\,
	\bigr\}\Big)~=~0 \eeq
	for all $k\geq 1$.
Here $C_0\doteq c\, 2^{-\frac{1}{\alpha}}\, \kappa_{1}^{1-\frac{1}{\alpha}}$. 
\end{lemma}
\v
{\bf Proof.}  {\bf 1.} Given a positive measure $\mu$ and a corresponding solution $u$ of 
(\ref{rm2}), consider the set 
$$A~\doteq~
	\bigl\{x\in \Omega_0\, ; \qquad  u(x)~\geq~C_0 |x|^{1/\alpha}\bigr\}.$$
Moreover, let 
\bel{phif}\tilde\mu~\doteq~\chi_A\cdot\mu
\eeq
be the measure obtained from $\mu$ by removing all the mass that lies outside $A$.

By (\ref{Idif}) it follows
\bel{Idif55}c\I^\alpha(\mu) - c\I^\alpha(\tilde \mu)~\geq~C_0 \int_{\Omega_0\setminus A} |x|^{1/\alpha}\,d\mu, \eeq
\v
{\bf 2.} 
To estimate the difference in the harvest functional, let $\tilde u$ be a solution 
to the same elliptic problem (\ref{rm2}) with $\mu$ replaced by $\tilde \mu$, such that 
\bel{sbsol}u(x)~\leq~\tilde u(x)~\leq~M\, \qquad\forall x\in \ov \Omega\,.\eeq
We compute
\bel{Hdif}\H(u, \mu ) - \H(\tilde u, \tilde \mu)~=~\int_{\Omega_0} u\, d\mu - \int_{\Omega_0} \tilde u\, d\tilde \mu~\leq~\int_{\Omega_0\setminus A} u\, d\mu~\leq~
C_0 \int_{\Omega_0\setminus A} |x|^{1/\alpha}\,d\mu\,.\eeq
Comparing  (\ref{Idif55}) with (\ref{Hdif}) we thus obtain
\bel{newmax} \H(u, \mu) - c\I^\alpha(\mu)~
\leq~\H(\tilde u, \tilde \mu) - c\I^\alpha(\tilde \mu).\eeq
Recalling that $u\leq \tilde u$, by (\ref{phif}) it follows
$$\tilde\mu\Big(\bigl\{x\in\Omega_0\, ; \quad \tilde u(x)~<~C_0 |x|^{1/\alpha}\,
	\bigr\}\Big)~\leq~
	\tilde\mu\Big(\bigl\{x\in \Omega_0\, ; \quad  u(x)~<~C_0 |x|^{1/\alpha}\,
	\bigr\}\Big)~=~0. $$
\v
{\bf 3.} If now $(u_k,  \mu_k)_{k\geq 1}$ is any maximizing sequence,
for every $k\geq 1$ we define
$$A_k~\doteq~
	\bigl\{x\in \ov\Omega\, ; \qquad  u_k(x)~\geq~C_0 |x|^{1/\alpha}\bigr\},\qquad
	\tilde\mu_k~\doteq~\chi_{A_k}\cdot\mu_k\,.$$
	Moreover, we let $\tilde u_k\geq u_k$ 
	be a solution to (\ref{rm2}) corresponding to the measure $\tilde \mu_k$.
	By the previous analysis,
 $(\tilde u_k, \tilde \mu_k)_{k\geq 1}$ is another maximizing sequence, satisfying 
(\ref{lbdu}).
\endproof
\v
When $f$ satisfies the assumptions (\ref{fp1}), any solution $u$ of (\ref{emu}) will take
values inside the interval $[0,M]$.
By the previous arguments, it thus  follows the existence of a maximizing sequence
$(u_k, \mu_k)_{k\geq 1}$, where the measures $\mu_k$  satisfy
\begi
\item[(i)]  Supp$(\mu_k)\subset\ov B_{r_0}$, where $C_0 r_0^{1/\alpha} = M$.
\item[(ii)] $\I^\alpha(\mu_k)\leq C$.
\endi
In particular, 
for every $r>0$ one has
\bel{mkb}
\sup_{k\geq 1} ~\mu_k\bigl(\{x\in \R^d\,;~|x|>r\}\bigr)~<~+\infty.\eeq
This does not necessarily imply that the total mass of the measures $\mu_k$ 
is uniformly bounded.
Indeed, they may concentrate more and more mass close to the origin.

To achieve the existence of an optimal measure, we thus need to work in the wider
class of positive measures $\mu$  on  the domain $\Omega_0$ in (\ref{om0}),
possibly with  infinite total mass.   As a preliminary, the definition of irrigation plan and irrigation cost
must  be extended to these more general measures.

If $\mu(\R^d)=+\infty$, an irrigation plan for $\mu$ is a map $\chi:\R_+\times\R_+
\mapsto\R^d$ with the properties (i)--(iii) introduced in Section~\ref{s:0}.
For every $m\geq 1$, call $\chi_m:[0,m]\times\R_+\mapsto\R^d$ 
the restriction of  $\chi$ to $[0,m]$.   Then the cost of $\chi$ is defined as
\bel{chic}
\E^\alpha(\chi)~\doteq~\lim_{m\to +\infty} \E^\alpha(\chi_m)
~=~\sup_{m\geq 1} 
\int_0^m\int_0^{\tau(\theta)} 
 \bigl|\chi_m(\theta,t)
\bigr|_{\chi_m}^{\alpha-1} \, dt
d\theta.\eeq
On the other hand, the harvest functional is defined as
\bel{hfc}
\H(u,\mu)~\doteq~\sup_{\ve\to 0+}~\int_{|x|> \ve}
u\,d\mu\,.\eeq
It is clear that the right hand sides  
in (\ref{chic}) and (\ref{hfc}) are well defined, possibly taking the value $+\infty$.
\v

We can now state our main result on the existence of an 
 optimal measure.
\v
\begin{theorem}\label{t:37}
Let the function $f$ satisfy  the assumptions (A1).  Then the maximization
problem for roots {\bf (OPR)} has an optimal solution $(u,\mu)$, where
$\mu$ is a  positive measure  on the domain $\Omega_0$ defined at (\ref{om0}).
The optimal measure $\mu$ has bounded support, but possibly unbounded total mass.
\end{theorem}

{\bf Proof.} 
{\bf 1.} Let $(u_k, \mu_k)$ be a maximizing sequence.  By the previous analysis we can assume that  all measures $\mu_k$ are supported inside a fixed ball $\ov B_\rho$, and  the quantities 
$\I^\alpha(\mu_k)$, $\H(u_k,\mu_k)$ remain uniformly bounded.

By possibly selecting a subsequence, we can assume the existence
of a positive measure $\mu$ such that the weak convergence 
$\mu_k\wto \mu$ holds on $\R^d\setminus\{0\}$. 
In other words,
$$\int \vp\, d\mu_k~\to~\int \vp\, d\mu$$
for every continuous function $\vp\in \C_c(\R^d\setminus\{0\})$.

Let $\mu^\ve$ be the restriction of $\mu$ to the subset $\{|x|>\ve\}$. Then
\bel{lia}\I^\alpha(\mu)~=~\sup_{\ve>0} ~\I^\alpha(\mu^\ve).\eeq
On the other hand, calling $\mu^\ve_k$ the restriction of $\mu_k$ 
to the set $\{ |x|\geq \ve\}$, the lower semicontinuity of the irrigation cost for
bounded measures implies 
\bel{lie}\I^\alpha(\mu^\ve)~\leq~\liminf_{m\to\infty}~ \I^\alpha(\mu^\ve_k)~\leq~
\liminf_{m\to\infty}~ \I^\alpha(\mu_k).\eeq
Combining (\ref{lia}) with (\ref{lie}) we obtain
\bel{li4}\I^\alpha(\mu)~\leq~\liminf_{m\to\infty} ~\I^\alpha(\mu_k).\eeq
\v
{\bf 2.} To complete the existence proof, we need to find a solution $u$ of (\ref{rm2})
and show  that 
\bel{husc}
\H(u, \mu)~\geq~\limsup_{k\to\infty}~\H(u_k, \mu_k).\eeq
Toward this goal,
choose radii $\rho_n\to 0$ such that 
\bel{mk0}\mu\Big(\{x\in\Omega_0\,;~~|x|=\rho_n\}\Big)~=~0\eeq
for all $n\geq 1$.   
Let $\mu^{\rho_n}_k, \mu^{\rho_n}$ be the restrictions of the measures $\mu_k,\mu$
to the closed sets $$V_n~\doteq ~\{x\in\Omega_0\,;~~|x|\geq\rho_n\}.$$
Thanks to (\ref{mk0}), for each $n\geq 1$ we have the weak convergence
$\mu^{\rho_n}_k\wto \mu^{\rho_n}$.
Let $u^{\rho_n}_k, u^{\rho_n}$ be the corresponding solutions of (\ref{rm2}).
By the analysis in \cite{BSu}, since all measures $\mu_k^{\rho_n}$ have uniformly bounded 
mass, for each $n\geq 1$ we have
\bel{uscu}
\int u^{\rho_n}\, d\mu^{\rho_n}~\geq~\limsup_{k\to \infty}\int u^{\rho_n}_k\, d\mu^{\rho_n}_k\,.
\eeq 
By Corollary~\ref{c:33} it follows
\bel{urn}\int u_k\, d\mu_k~\leq~\int u^{\rho_n}_k\, d\mu^{\rho_n}_k+\eta(\rho_n).\eeq
\v
{\bf 3.} 
We now observe that, as $n\to\infty$, the sequence of measures $\mu^{\rho_n}$
is increasing while the sequence of solutions $u^{\rho_n}$
is  decreasing.   Setting 
$$u(x)~\doteq~\inf_{n\geq 1}u^{\rho_n}(x),$$
one checks that $(u,\mu)$ is a solution to (\ref{rm2}).
Moreover, for any $\delta>0$ one has
\bel{udmu}\int_{|x|\geq\delta} u\, d\mu~=~\inf_{n\geq 1} \int_{|x|\geq\delta} u^{\rho_n}
\, d\mu
.\eeq
Given $\ve>0$, we can find $\delta>0$ and then an integer $\Hat n$ large enough  so that 
\bel{eds}\rho_{\Hat n}~\leq~\delta,\qquad 
\eta(\rho_{\Hat n})~\leq~ \eta(\delta)~<~\ve,\qquad\qquad 
\int_{|x|\geq\delta} u^{\rho_{\Hat n}}
\, d\mu~<~\int_{|x|\geq\delta} u\, d\mu+\ve\,.\eeq
  Using (\ref{urn}), then  (\ref{uscu}), and finally
(\ref{eds}),  we conclude
$$\bega{l}\ds\limsup_{k\to\infty} \int u_k \, d\mu_k~\leq~\limsup_{k\to\infty}
\int  u_k^{\rho_{\Hat n}} 
\, d\mu_k^{\rho_{\Hat n}}  +\eta(\rho_{\Hat n})
~\leq~\int  u^{\rho_{\Hat n}} 
\, d\mu^{\rho_{\Hat n}}  +\eta(\rho_{\Hat n})\\[4mm]
\ds\qquad < ~\left(\int_{|x|\geq\delta}u\, d\mu +\ve\right) +  \ve~\leq~\int u\, d\mu +2\ve.
\enda$$
Since $\ve>0$ was arbitrary, this completes the proof.
\endproof

\v

\section{Dependence on parameters}
\label{s:7}
\setcounter{equation}{0}
Let $\alpha\in [0,1]$ be given.   In Sections 3 and 4
we have proved the existence of an optimal configuration of 
tree roots and tree branches, where the optimal measure $\mu$ has bounded support.  
Here we are interested in how  this support
depends on parameters.
More precisely, given a measure $\mu$ on $\R^d$,
let 
\bel{RSUP}
R(\mu)~\doteq~\inf\Big\{ r>0\,;~~\mu\bigl(\{|x|>r\}\bigr)=0\Big\}\eeq
be the radius of the smallest ball centered at the origin which contains
the support of $\mu$. 

We first consider the optimization problem {\bf (OPB)} for tree branches.
We seek an upper bound on $R(\mu)$, 
depending on the dimension $d$, the  constants $\alpha,c$, and the
$\L^1$ norm of the function $\eta$  in (\ref{OB}), measuring the intensity of light from various directions.   

As a preliminary, we recall  how the irrigation cost behaves under rescalings.
Given a measure $\mu$ and a constant
 $\lambda>0$, we  define the measures $\lambda\mu$ and $\mu^\lambda$ 
respectively by setting
\bel{RSD}(\lambda \mu)(A)~\doteq~\lambda\, \mu(A),\qquad\qquad \mu^\lambda(A)~\doteq~\mu(\lambda^{-1}A),\eeq
for every Borel set $A\subseteq\R^d$.

\v
\begin{lemma}\label{RS0} For any positive Radon measure $\mu$ on $\R^d$ and any $\lambda>0$, 
$0\leq\alpha\leq 1$, the following holds:
\bel{IR} \I^{\alpha}(\lambda\mu)~=~\lambda^{\alpha}\I^{\alpha}(\mu),\qquad\qquad 
\I^{\alpha}(\mu^{\lambda})~=~\lambda\I^{\alpha}(\mu).\eeq
\end{lemma}
\v

{\bf Proof.} {\bf 1.} To prove the first identity in (\ref{IR}), let $\Theta=[0,M]$ and let $\chi:[0,M]
\times\R_+\mapsto\R^d$ be an admissible irrigation plan for $\mu$.
Then the map $\chi^\lambda:[0,\lambda M]
\times\R_+\mapsto\R^d$, defined by 
\bel{hpg}\chi^\lambda( \theta, t)  ~=~\chi(\lambda^{-1}\theta, t)
\eeq
is an admissible irrigation plan for $\lambda\mu$.  Its cost is computed by
\bel{change}
\bega{rl}\ds \E^\alpha(\chi^\lambda)&\ds=~\int_{[0, \lambda M]}\Big(\int_{\R_{+}}|\chi^\lambda(\theta,t)|_{\chi^\lambda}^{\alpha-1}\cdot|{\dot\chi}^\lambda(\theta,t)|\, 
dt\Big)\, d\theta\\[4mm]
\ds &\ds=~\int_{[0,M]}\Big(\int_{\R_{+}}\Big(\lambda |\chi(\theta,t)|_{\chi}\Big)^{\alpha-1}\cdot|\dot{\chi}(\theta,t)|\, dt\Big)\,\lambda\, d\theta\\[4mm]
\ds &\ds=~\lambda^\alpha\E^\alpha(\chi).
\enda\eeq
Taking the infimum over all irrigation plans we obtain
 $\I^\alpha(\lambda\mu)\leq\lambda^\alpha\I^\alpha(\mu)$. 
 Replacing $\lambda$ by $\lambda^{-1}$ we obtain the opposite inequality.
\v
{\bf 2.} To prove the second identity, consider any $\lambda>0$ and 
let  $\chi:\Theta\times\R_+\mapsto\R^d$ be an  
irrigation plan for $\mu$.  Then 
$\chi^\dagger:\Theta\times\R_+\mapsto\R^d$
defined by
$$\chi^\dagger(\theta,t) ~=~\lambda \cdot\chi(\lambda^{-1}\theta,t)$$
is an admissible irrigation plan for $\mu^\lambda$. 
Performing the change of variables $\tilde\theta=\lambda^{-1}\theta$,
 its cost is computed as
\bel{change2}
\bega{rl}\ds \E^\alpha(\chi^\dagger)&\ds=~\int_\Theta\Big(\int_{\R_{+}}|\chi^\dagger(\theta,t)|_{\chi^\dagger}^{\alpha-1}\cdot|{\dot\chi}^\dagger(\theta,t)|\, 
dt\Big)\, d\theta\\[4mm]
\ds &\ds=~\int_\Theta\Big(\int_{\R_{+}}|\chi(\tilde \theta,t)|_{\chi}^{\alpha-1}\cdot
 |\dot{\chi}(\tilde \theta,t)|\, dt\Big)\,\lambda\, d\tilde \theta\\[4mm]
\ds &\ds=~\lambda\E^\alpha(\chi).
\enda\eeq
Taking the infimum over all irrigation plans we obtain
 $\I^\alpha(\mu^\lambda)\leq\lambda\I^\alpha(\mu)$. 
 Replacing $\lambda$ by $\lambda^{-1}$ we obtain the opposite inequality.\endproof

\v
Similar formulas relate the sunlight captured by a rescaled measure.
Namely, as proved in \cite{BSu}, one has
\bel{SR}
\S^{b\eta}(\mu)~=~b\,\S^\eta(\mu),
\qquad\qquad \S^\eta(\lambda^{d-1} \mu^\lambda) ~=~\lambda^{d-1} \S^\eta(\mu).\eeq
\v
Thanks to the rescaling properties (\ref{IR}) and (\ref{SR}), 
the solution to the problem
\bel{P1}\hbox{maximize:}\quad \S^\eta(\mu)-\I^\alpha(\mu)\eeq
can be related to the solutions to the family of problems
\bel{P2}\hbox{maximize:}\quad \S^{b\eta}(\mu)-c\I^\alpha(\mu),\eeq
for any constants $b,c>0$.
\v
\begin{lemma}\label{p:62}
Assume $1 - \frac{1}{d - 1}\doteq \alpha^*< \alpha$ and assume that the measure $\mu$ is optimal for the problem 
(\ref{P1}).  Then, for any given constants $b, c>0$, the measure 
\bel{sclmu}\tilde \mu ~=~ \lambda^{d - 1} \mu^{\lambda}\, , \qquad \lambda~\doteq~\Big(\frac{b}{c}\Big)^{\frac{1}{1 + (\alpha - 1)(d - 1)}}\, \eeq
provides an optimal solution to the problem (\ref{P2}).
\end{lemma}

{\bf Proof.} Given any measure $\mu$, define $\tilde \mu$ by setting
\bel{RS1}\lambda^{d - 1}\mu^{\lambda}~=~\tilde \mu\, .\eeq
By the rescaling formulas (\ref{IR}) and (\ref{SR}), one has
\bel{RS2}\bega{rl}\ds \S^{b\eta}(\tilde \mu) - c\I^\alpha(\tilde \mu)&\ds =~\S^{b\eta}(\lambda^{d-1} \mu^\lambda) - c\I^\alpha( \lambda^{d - 1}\mu^{\lambda})\\[4mm]
&\ds =~b\lambda^{d-1}\S^\eta(\mu) - c\lambda^{1 + \alpha(d - 1)}\I^\alpha(\mu).
\enda\eeq
By the definition of $\lambda$ in (\ref{sclmu}), it follows
\bel{RS3}\S^{b\eta}(\tilde \mu) - c\I^\alpha(\tilde \mu)~=~b^{\frac{1 + \alpha(d -1)}{1 + (\alpha - 1)( d -1)}} c^{\frac{1-d}{1 + (\alpha -1 )( d -1)}}\Big( \S^\eta(\mu) - \I^\alpha(\mu)\Big).\eeq
Therefore, $\S^{b\eta}(\tilde \mu) - c\I^\alpha(\tilde \mu)$ attains the maximum possible value 
if and only if $ \S^\eta(\mu) - \I^\alpha(\mu)$ attains the  maximum possible value. This completes our proof.
\endproof
\v
Our next result provides an estimate on the size of the support of the optimal measure $\mu$.
\v
\begin{proposition} In the same setting as  Theorem~\ref{t:21}, for any $d\geq 2$ and $1-{1\over d-1}<\alpha<1$, there is a 
constant $C_{\alpha,d}$ such that any optimal measure $\mu$ for the problem (\ref{OB}) is supported inside
a ball of radius
\bel{RM}R(\mu)~\leq~{C}_{d,\alpha} \Big(\frac{\|\eta\|_{\mathbf{L}^1}}{c}\Big)^{\frac{1 }{1 + (\alpha - 1)(d - 1)}} \, .\eeq
When $\alpha=1$ one simply has
\bel{dop7}R(\mu)~\leq~\frac{\|\eta\|_{\mathbf{L}^1}}{c}\, .\eeq
\end{proposition}

{\bf Proof.}   {\bf 1.}   Consider first the special  case where $\|\eta\|_{\L^1} = c = 1$.
By (\ref{sf1})-(\ref{ria}) and (\ref{Ibb}), we then have
\bel{dop1}\Tilde C_{d,\alpha}\Big(\I^\alpha(\mu)\Big)^{\frac{d-1}{1+\alpha(d-1)}}- \I^\alpha(\mu)~\geq~0,\eeq
where $\Tilde C_{d,\alpha}$ is a constant which only depends on $d$ and $\alpha$. 
Therefore, in (\ref{k*}) one can take  the constant
\bel{dop2}\kappa_1~ \doteq ~\Tilde C_{d,\alpha}^{~\frac{1 + \alpha( d - 1)}{1 + (\alpha - 1)(d - 1)}}.\eeq
\v
{\bf 2.} In the case  $1 - \frac{1}{d-1}<\alpha < 1$, 
we choose the radius
\bel{dop4}r_1~\doteq~\alpha^{\alpha\over\alpha-1}\,\kappa_1\,.\eeq
By (\ref{uvw0})-(\ref{uvw11}), this yields
\bel{dop3} |x|\geq r_1\qquad\implies\qquad\alpha |x|_{\chi}^{\alpha -1 }~\geq~\alpha \Big(\frac{\kappa_1}{r_1}\Big)^{1 - {1\over \alpha}}~\geq~1\, .\eeq
By the argument following (\ref{oenway}), the optimal measure is supported on the 
ball $\ov B_{r_2}$, where the radius $r_2$ satisfies 
\bel{dop5}r_2 - r_1~=~1\,.\eeq
Recalling 
 (\ref{dop4}) and (\ref{dop2}), we obtain an upper bound on $R(\mu)$, namely
\bel{dop6} R(\mu)~\leq\ds ~r_2~=~r_1+1 ~=~
\alpha^{\alpha\over\alpha-1}\,\kappa_1+1
~=~\ds \alpha^{\alpha\over\alpha-1}\,\Tilde C_{d,\alpha}^{~\frac{1 + \alpha( d - 1)}{1 + (\alpha - 1)(d - 1)}}+1 ~\doteq~C_{d,\alpha}\,.\eeq
\v
{\bf 3.}  To cover  the general case, let $b\doteq \|\eta\|_{\L^1}$.   
Then we can write $\eta = b\, \tilde \eta$, where $\|\tilde\eta\|_{\L^1}=1$.  
By Proposition~\ref{p:62},  a measure $\tilde \mu$ is optimal for the problem
$$\hbox{maximize:}\qquad \S^{\tilde\eta}(\mu) - \I^\alpha(\mu)$$
if and only if the measure
\bel{muti}\mu ~\doteq~\lambda^{d-1} \tilde\mu^\lambda,\qquad\qquad \lambda~=~\left( b\over c\right)^{1\over 1+(\alpha-1)(d-1)}
\eeq
is optimal for the problem (\ref{OB}).   

Since $\|\tilde\eta\|_{\L^1} =1$, by the previous step the measure $\tilde \mu$ is supported on a ball of radius
$\R(\tilde\mu)\leq C_{d,\alpha}$.   In turn, by (\ref{muti}), the measure $\mu$ is supported on a ball of radius
$R(\mu)=\lambda R(\tilde\mu)\leq\lambda\, C_{d,\alpha}$.    This proves (\ref{RM}).
\v
{\bf 4.} When $\alpha = 1$, the estimate (\ref{dop7}) is an immediate consequence of   (\ref{r2}).
\endproof
\v
\begin{remark} {\rm The radius of the smallest ball containing the support of $\mu$ 
can be regarded as the ``size" of the tree.   
As expected, the above analysis indicates that the optimal size increases with the 
amount of sunlight $\|\eta\|_{\L^1}$, and decreases with the factor $c$ multiplying the 
irrigation cost.}
\end{remark}

Similar questions can be asked in connection with the optimization problem {\bf (OPR)} for tree roots.  
More precisely,
assume that the
diffusion depends on a parameter $\sigma>0$,
and let a function $f:\R_+\mapsto\R$ be given, as in (\ref{fp1}).  
Consider the optimization problem
\bel{m5}
\hbox{maximize:}\quad \H(u,\mu)-c\I^\alpha(\mu),\eeq
\bel{m6}\hbox{subject to:}\quad 
\left\{ \bega{rl}\sigma\, \Delta u+ a\,  f(bu)-u\mu&=~0,\qquad x\in \R^d_-\,,\\[3mm]
u_{x_d}&=~0,\qquad x_d=0.\enda\right.\eeq

Let $\mu$ be an optimal measure and call $R(\mu)$ the radius of the smallest 
closed ball, centered at 
the origin, which contains the support of $\mu$.
We wish to understand how this radius depends on the parameters $a,b,c$, and $\sigma$.

Throughout the following, we assume that $d \geq 2$ and $1 - \frac{1}{d} < \alpha \leq 1$, while
 $f$ satisfies (A1).

As a first step, we consider the problem 
\bel{m7}
\hbox{maximize:}\quad \H(\tilde u,\tilde \mu)-\tilde c\,\I^\alpha(\tilde \mu),\eeq
\bel{m8}\hbox{subject to:}\quad 
\left\{ \bega{rl} \Delta\tilde  u+  f(\tilde u)-\tilde u\,\tilde \mu&=~0,\qquad x\in \R^d_-\,,\\[3mm]
\tilde u_{x_d}&=~0,\qquad x_d=0,\enda\right.\eeq
and prove a rescaling result, similar to Proposition~\ref{p:62}.

\v
\begin{lemma}\label{p:64}	A couple $(\tilde u, \tilde \mu)$ is an optimal solution 
to (\ref{m7})-(\ref{m8}) if and only if  $(u, \mu)$ 
is an optimal solution to (\ref{m5})-(\ref{m6}), where
	\bel{newop1} u(x) \,\doteq\,\frac{1}{b}\tilde u(\lambda^{-1}x)\, , 
	\quad \qquad  \mu\,\doteq\,\sigma\lambda^{d - 2}\tilde \mu^{\lambda}\,,\qquad\quad \lambda\doteq\sqrt{\frac{\sigma}{ab}}\,, \qquad \tilde c ~\doteq~\frac{c\sigma^\alpha }{a } \lambda^{ 1-(1- \alpha )d  - 2\alpha}.\eeq
\end{lemma}

{\bf Proof.} {\bf 1.} Let $(\tilde u,\tilde \mu)$ be an optimal solution to (\ref{m8}).
For any test function $\Tilde \vp\in \C_c(\R^d)$, set $ \vp(x)\doteq \Tilde \vp(\lambda^{-1} x)$.
By the definition of the 
rescaled measure $\tilde \mu^\lambda$ in (\ref{RSD}),  one has  
$$\int \vp\, d\tilde \mu^\lambda ~=~ \int \Tilde 
\vp\, d\tilde \mu\,.$$
Therefore
\bel{E1}\bega{c}\ds \int_{\R^d_-} \sigma\, \nabla u(x)\, \nabla \varphi(x)\, dx ~\ds =~\frac{\sigma}{b\lambda}\int_{\R^d_-}  \nabla\tilde  u(\lambda^{-1} x)\, \nabla\Tilde \varphi(\lambda^{-1}x)\, dx~=~\frac{\sigma \lambda^{d -2 }}{b}\int_{ {\R^d_-}} \nabla 
\tilde u(y)\cdot \nabla \Tilde \varphi(y)\, dy\, ,
\\[4mm]
\ds \int_{\R^d_-} a\cdot f(b u(x))\, \varphi(x)\, dx~\ds =~ \int_{\R^d_-} a\,f(u(\lambda^{-1} x))\, 
\Tilde \varphi(\lambda^{-1}x)\, dx ~=~a\lambda^d\int_{\R^d_-}  f(u(y))\,\Tilde \varphi(y)\, dy\, ,
\\[4mm]
\ds \int_{\ov\R^d_-} u(x)\, \varphi(x)\, d\mu~ =~
\int_{\ov \R^d_-}\frac{\sigma\lambda^{d -2}}{b}\tilde u(\lambda^{-1} x)\, \Tilde \varphi(\lambda^{-1}x)\, 
d\tilde \mu^{\lambda}~=~\frac{\sigma\lambda^{ d -2}}{b}\int_{ \ov \R^d_-}\tilde  u(y)\, \Tilde \varphi(y)\, d\tilde \mu\, .
\enda\eeq
Since $(\tilde u,\tilde \mu)$ is a solution to (\ref{m8}), we have
\bel{E2}
-\int_{ {\R^d_-}} \nabla \tilde u(y)\cdot \nabla \Tilde \varphi( y)\, dx + \int_{\R^d_-}  f(\tilde u(y))\,
\Tilde \varphi(y)\, dx-\int_{ {\R^d_-}} \tilde u(y)\,\Tilde \varphi(y)\, d\mu~=~0\,.
\eeq
Taking  $\lambda\doteq \sqrt{\frac{\sigma}{ab}}$, from (\ref{E1})-(\ref{E2}) it follows
\bel{E3}- \sigma \int_{\R^d_-} \nabla u(x)\cdot\nabla\varphi(x)\, dx + \int_{\R^d_-} a\cdot f(b u(x))\, \varphi(x)\, dx-\int_{\ov \R^d_-} u(x)\, \varphi(x)\, d\mu~=~0\, .\eeq
Since the test function $\varphi$ is arbitrary, we conclude that $( u,  \mu)$ is a solution to (\ref{m6}).
\v
{\bf 2.} We now claim that $(u,  \mu)$ is actually a solution to the optimization problem 
(\ref{m5})-(\ref{m6}). Indeed, given any measure $ \nu$, there is a unique measure $\tilde \nu$ such that
\bel{E4}{\nu}~=~\sigma\lambda^{d-2 }\tilde\nu^{\lambda}\, ,\qquad \lambda~\doteq~\sqrt{\frac{\sigma}{ab}}\, .\eeq
By the preceding argument, if $\tilde v$ is a solution to (\ref{m8}) corresponding to the measure $\tilde\nu$, 
then $v(x) \doteq b^{-1}\tilde v(\lambda^{-1} x)$ is a solution to (\ref{m6}) corresponding to the measure 
$ \nu$.  Moreover,
\bel{E5}
\H( v, \nu)~=~\int_{ {\R^d_-}} a\cdot f(b v(x))\, dx~=~a\int_{ {\R^d_-}}f(\tilde v(\lambda^{-1} x
))\, dx~=~a\lambda^d\, \H(\tilde v,\tilde \nu)\, .
\eeq
On the other hand, by  (\ref{IR}) it follows
\bel{E6} c\I^\alpha ( \nu)~=~c\I^{\alpha}(\sigma \lambda^{d -2}\tilde \nu^\lambda)~=~c\sigma^{\alpha}
\lambda^{\alpha(d - 2) + 1}\, \I^\alpha (\tilde \nu)\, .\eeq
By the choice of $\tilde c$ in (\ref{newop1}), from  (\ref{E5})-(\ref{E6}) we conclude
\bel{E7}\H( v, \nu) - c\I^\alpha( \nu)~=~a\lambda^d \Big( \H(\tilde v, \tilde \nu) - \tilde c \I^\alpha(\tilde \nu)\Big).\eeq
Therefore, $  \H(v, \nu) - c\I^\alpha(\nu)$ attains its maximum if and only if  
$\H(\tilde v,\tilde \nu) - \tilde c\I^\alpha(\tilde \nu) $ attains its maximum. This completes the proof.
\endproof

\v

\begin{lemma}\label{S:65}
	Let $(u, \mu)$ be an optimal solution for (\ref{OR})-(\ref{rm2}). 
Then there is a constant $C$, depending  on 
$d$, $\alpha$ and $f$, such that $\mu$ is supported inside the ball of radius
	\bel{ps1}R(\mu)~\leq~C\cdot c^{\frac{- 1}{ 1 - ( 1 - \alpha)d}}\, .\eeq
	In the special case where $\alpha = 1$, one has the simpler estimate
	\bel{ps0} R(\mu)~\leq~\frac{M}{c}\,.\eeq
\end{lemma}
\v

{\bf Proof. 1. } Assume  that $1 -\frac{1}{d}< \alpha < 1$. By (\ref{aasse}), there is a constant $C_0$, 
depending 
on $d$, $\alpha$, and $f$, such that
$$\H(u, \mu) - c\I^\alpha (\mu)~\leq~C_0\Big(\I^\alpha (\mu)\Big)^{\frac{d}{ 1 + \alpha d}} - c\I^\alpha(\mu)\, .$$
Requiring $C_0\Big(\I^\alpha (\mu)\Big)^{\frac{d}{ 1 + \alpha d}} - c\I^\alpha(\mu)\geq 0$, 
one obtains an a priori bound $\kappa$ for the irrigation cost $\I^\alpha$, namely
\bel{kappabd} \I^\alpha(\mu)~\leq~\kappa~\doteq~C_1\cdot c^{- \frac{1 + \alpha d}{ 1 + (\alpha - 1)d}}\, ,\eeq
for some constant $C_1$ depending on $d$, $\alpha$ and $f$.

For any $\gamma > 0 $, by same argument as in (\ref{uvw0}) we can find a radius $r_1$ such that
\bel{such} |x|~\geq~r_1\, ,\qquad\Longrightarrow\qquad \alpha |x|_\chi^{\alpha -1}~\geq~\gamma\, .\eeq
 Indeed, by (\ref{uvw0}) one can choose
\bel{such1}r_1~\doteq~\alpha^{\frac{\alpha}{\alpha - 1}}\, \gamma^{\frac{\alpha}{1 -\alpha}}\, \kappa\, .\eeq 

We now split $\mu = \mu^\flat + \mu^\sharp$ as in (\ref{md11}), choosing $r_2$ so that
\bel{such2}c\gamma (r_2 - r_1) ~=~M\, , \qquad r_2~=~r_1 + \frac{M }{c\gamma}\,.\eeq
Using (\ref{oenway}) and (\ref{such}), 
 we now obtain
$$\I^\alpha(\mu^\flat + \mu^\sharp)~\geq~\I^\alpha(\mu^\flat) + \gamma\,\mu^\sharp(\R^d).$$
By (\ref{kappabd}), (\ref{such1}), and (\ref{such2})  it follows
\bel{such3} r_1~=~\frac{C_2\cdot \gamma^{\frac{\alpha}{1 -\alpha}}}{c^{\frac{1 +\alpha d}{1 + (\alpha -1)d}}}\,,
\qquad\qquad
r_2~\doteq~ r_1 + \frac{M}{c\gamma}~=~\frac{C_2\cdot \gamma^{\frac{\alpha}{1 -\alpha}}}{c^{\frac{1 +\alpha d}{1 + (\alpha -1)d}}} + \frac{M}{c\gamma}\, .
\eeq
To achieve the best estimate on the radius $r_2$, we minimize the right hand side of 
(\ref{such3}) over all possible choices of $\gamma> 0 $. 
Taking $\gamma\doteq \Big(\frac{M}{C_2}\Big)^{1 -\alpha} c^{\frac{(1 -\alpha) d}{1 + (\alpha - 1)d}}$, one
obtains
\bel{such5}r_2~=~2\frac{M^\alpha\cdot C_2^{1 -\alpha}}{c^{\frac{1}{1 + (\alpha - 1)d}}}~=~C\cdot c^{\frac{- 1}{1 - (1 -\alpha)d}}\, ,\eeq
where the constant $C$ only depends on  $d$, $\alpha$ and $f$.
\v
{\bf 2.} Next, if $\alpha = 1$, by (\ref{mdk}) one has
\bel{ps4}\I(\mu^\flat + \mu^\sharp) - c\Big(\I(\mu^\flat)\Big)~=~c\I(\mu^\sharp)~\geq~cr\mu^\sharp(\R)\, .\eeq
On the other hand, by the assumptions in (A1),   any solution $u$ of (\ref{emu}) takes values inside $[0,M]$.
Therefore, any measure containing some mass outside the ball $\ov B_\rho$, centered at the origin
with radius $\rho = M/c$, cannot be optimal. \endproof
\v
Combining Lemmas~\ref{p:64} and \ref{S:65}, we now obtain an estimate on the support of the 
optimal measure for the general optimization problem for tree roots. 

\begin{proposition} \label{P57}
Assume that $d \geq 2$ and $1 - \frac{1}{d} < \alpha \leq 1$, 
	and let $f$ satisfy the assumption $(A1)$.   Then
there is a constant $C$ only depending on the dimension $d$, $\alpha$ and $f$ such that any optimal measure $\mu$ for the problem (\ref{m5})-(\ref{m6}) is supported inside a ball of radius
	\bel{ps2}R(\mu)~\leq~C \,a^{\frac{1 -\alpha}{ 1 - (1-\alpha )d}}\,\, b^{\frac{- \alpha}{ 1 - (1-\alpha )d}}\, c^{\frac{- 1}{ 1 - (1-\alpha )d}}\,. \eeq
	When $\alpha = 1$ one simply has
	\bel{ps3}
	R(\mu)~\leq~{M\over bc}\,.
	\eeq
\end{proposition}

{\bf Proof. } {\bf 1.} Consider first the case $\alpha < 1$. 
Let $(u,\mu)$ be an optimal solution to (\ref{m5})-(\ref{m6}).
Then by Lemma~\ref{p:64} the couple $(\tilde u, \tilde \mu)$ in (\ref{newop1})
provides an optimal  solution to (\ref{m7})-(\ref{m8}).

Using  Lemma~\ref{S:65} and performing the variable transformations in  (\ref{newop1}), this yields 
\bel{ps8} R(\mu)~=~\lambda R(\tilde\mu)~\leq~
\sqrt{ \sigma\over ab} \cdot C\,\tilde c ^{\frac{-1}{1 + (\alpha - 1) d}}
~=~C\, \sqrt{ \sigma\over ab}\cdot \left[ {c\sigma^\alpha\over a} 
\Big( {\sigma\over ab}\Big)^ {1-(1-a)d -2\alpha\over 2}\right] ^{-1\over 1-(1-\alpha) d}.\eeq
After some simplifications, from (\ref{ps8}) one obtains precisely   (\ref{ps2}).
\v
{\bf 2.} Since every solution $u$ of (\ref{m6}) satisfies $0\leq u(x)\leq {M\over b}$, the estimate (\ref{ps3}) is clear.
\endproof

\begin{remark} {\rm By assumption, in (\ref{ps2}) all denominators are positive:
$1+(\alpha-1)d>0$.  From Proposition~\ref{P57} it follows that the support 
of $\mu$ decreases as the factor  $c$ multiplying the transportation cost 
becomes larger.  Somewhat surprisingly, the diffusion coefficient $\sigma$ 
does not seem to play a major role in determining the optimal size of tree roots.
Indeed, on the right hand side of  (\ref{ps8}) the various powers of $\sigma$ exactly cancel each other.
}
\end{remark}
\v


\begin{thebibliography}{6111}


\v
\bibitem{BCM} M.~Bernot,  V.~Caselles, and J.~M.~Morel,
{\it  Optimal transportation networks. Models and theory.} 
Springer Lecture Notes in Mathematics {\bf 1955},
Berlin, 2009.

\bibitem{BCM1} M.~Bernot,  V.~Caselles, and J.~M.~Morel,
 The structure of branched transportation networks.
 {\it Calculus of Variations}  (2008), 279-317. 
 
\bibitem{BG1}
 L.~Boccardo and  T.~Gallou\"et,
 Non-linear elliptic and parabolic equations involving measure data.
{\it J.~Functional Analysis} {\bf  87} (1989), 149--169.
 
 \bibitem{BGO}
 L.~Boccardo and  T.~Gallou\"et, and L.~Orsina,
 Existence and uniqueness of entropy solutions for nonlinear elliptic equations with measure data.
{\it Ann. Institut H.~Poincar\'e Nonlin. Anal.} {\bf 13} (1996), 539--551.
 
\bibitem{BS} L.~Brasco and F.~Santambrogio, 
An equivalent path functional formulation of branched transportation problems. 
{\it Discrete Contin. Dyn. Syst.} {\bf 29} (2011), 845--871.

\bibitem{BCS} A.~Bressan,
G.~Coclite and W.~Shen,
A multi-dimensional optimal harvesting problem with measure valued solutions,
{\it SIAM J. Control Optim.} {\bf 51} (2013), 1186--1202.


\bibitem{BSu}  A.~Bressan,  and Q.~Sun,
On the optimal shape of tree roots and branches, submitted.
Available on
http://arxiv.org/abs/1803.01042

\bibitem{BDM} G.~Buttazzo and G.~Dal Maso,  Shape optimization for Dirichlet problems: relaxed formulation and optimality conditions.
{\it Appl. Math. Optim. } {\bf 23}  (1991),  17--49.


\bibitem{DMOP}
G.~Dal Maso, F.~Murat, L.~Orsina, and A.~Prignet,  Renormalized solutions of elliptic equations with general measure data.  {\it Ann. Scuola Norm. Sup. Pisa Cl. Sci.} {\bf   28}
  (1999),  741--808.





\bibitem{FZ} H.~Federer and W.~Ziemer,
The Lebesgue set of a function whose distribution derivatives are p-th power summable.
{\it  Indiana Univ. Math. J.} {\bf 22} (1972), 139--158.

\bibitem{EG}  L.~C.~Evans and R.~F.~Gariepy, 
{\it Measure Theory and Fine Properties of Functions}. CRC Press, 1992.

\bibitem{MMS} F.~Maddalena, J.~M.~Morel, and S.~Solimini, 
A variational model of irrigation patterns, 
{\it Interfaces Free Bound.} {\bf 5}  (2003), 391--415.

\bibitem{MS} F.~Maddalena and S.~Solimini,
Synchronic and asynchronic descriptions of irrigation problems. 
{\it Adv. Nonlinear Stud.} {\bf 13} (2013), 583--623. 

\bibitem{MoS} J.~M.~Morel and F.~Santambrogio, 
The regularity of optimal irrigation patterns. 
{\it Arch. Ration. Mech. Anal.} {\bf  195}
 (2010), 499--531.

 \bibitem{OS}
 E.~Oudet and F.~Santambrogio,
 A Modica-Mortola approximation for branched transport and applications. 
{\it  Arch. Rational Mech.  Anal.} {\bf  201} (2011), 115--142.



\bibitem{PSX} P.~Pegon, F.~Santambrogio, and Q.~Xia, 
A fractal shape optimization problem in branched transport. 
{\it J. Math. Pures Appl.}, to appear.

\bibitem{S1}
F.~Santambrogio, 
Optimal channel networks, landscape function and branched transport. 
{\it Interfaces Free Bound.} {\bf 9} (2007),  149--169. 

\bibitem{S2}  F.~Santambrogio,  A Modica-Mortola approximation for branched transport.  
 {\it C. R. Acad. Sci. Paris, Ser. I}, {\bf  348} (2010) 941--945.

\bibitem{X03}
Q.~Xia, Optimal paths related to transport problems, 
{\it Comm. Contemp. Math.} {\bf 5}  (2003), 251--279.

\bibitem{X4} Q.~Xia, 
Motivations, ideas and applications of ramified optimal transportation.
{\it ESAIM Math. Model. Numer. Anal.}
{\bf  49} (2015),  1791--1832. 

\end{thebibliography}
\end{document}